\documentclass[journal]{IEEEtran}
\usepackage{cite}
\usepackage{url}
\usepackage{empheq}
\usepackage{float}
\usepackage{dsfont}
\usepackage[hidelinks]{hyperref}

\usepackage[normalem]{ulem}
\usepackage{amsmath,amssymb,amsfonts}
\allowdisplaybreaks[4]
\usepackage{graphicx}
\usepackage{textcomp}
\usepackage{amsmath}
\usepackage{enumerate}
\usepackage{epstopdf}
\usepackage{array}
\usepackage{booktabs}
\usepackage{subfigure}
\usepackage{multirow}
\usepackage{soul, color}
\usepackage[usenames,dvipsnames]{xcolor}
\usepackage[latin1]{inputenc}
\usepackage{cases}

\usepackage{amsthm}
\newtheorem{theorem}{Theorem}

\newtheorem{proposition}{Proposition}

\newtheorem{definition}{Definition}

\soulregister\cite7 
\soulregister\citep7 
\soulregister\citet7 
\soulregister\ref7 
\soulregister\pageref7

\usepackage{bbding}

\usepackage{nomencl}
\makenomenclature
\usepackage{etoolbox}
\renewcommand\nomgroup[1]{%
  \item[\bfseries
  \ifstrequal{#1}{A}{Acronyms}{%
  \ifstrequal{#1}{S}{Symbols}{%
  \ifstrequal{#1}{U}{Units}{}}}%
]}

\usepackage{amsmath}
\allowdisplaybreaks[4]

\usepackage{algorithm}
\usepackage{algpseudocode}
\usepackage{bm}

\usepackage{bbm}

\begin{document}

\title{Online Carbon-Aware Power Flexibility Quantification of Distributed Data Centers}

\author{Xiaowei Wang, Yue Chen, \IEEEmembership{Senior Member, IEEE}
\thanks{X. Wang and Y. Chen are with the Department of Mechanical and Automation Engineering, The Chinese University of Hong Kong, Hong Kong, China. (E-mail:  xwwang@mae.cuhk.edu.hk, yuechen@mae.cuhk.edu.hk)}
}
\maketitle

\begin{abstract}
Geographically distributed data centers (DCs) can provide spatio-temporal power flexibility to electric grids through the coordinated management of computing resources and energy. Beyond their energy consumption, the carbon footprints of DCs are non-negligible, which constrain the power flexibility that DCs can offer to the power grid. To address this challenge, this paper proposes a carbon-aware power flexibility quantification framework for DCs. We first develop an optimization model that characterizes DC power flexibility in a grid-usable form, namely the Cartesian product of adjustable power demand intervals across time slots. An online algorithm based on Lyapunov optimization is then designed to solve the flexibility quantification model in real time under uncertain workload arrivals and fluctuating grid carbon intensity. Furthermore, an integrated disaggregation and feedback mechanism is devised to translate the power system operator's dispatch commands into real-time DC operating decisions and to capture their impact on future DC flexibility. Several desirable properties of the proposed model and method are proven, including disaggregation feasibility, bounded workload processing delays, and a bounded suboptimality gap. Numerical simulations demonstrate that our approach effectively unlocks DC power flexibility while satisfying carbon emission constraints, outperforming baseline algorithms.
%Geographically distributed data centers (DCs) can provide spatio-temporal power flexibility to electric grids through coordinated computing resources and energy management. However, their growing energy consumption has also made them major contributors to global carbon emissions. Consequently, characterizing DC power flexibility requires strictly accounting for carbon emission limits. This paper proposes a novel carbon-aware power flexibility quantification framework for DCs. We represent DC power flexibility through dynamic flexibility intervals and develop an online algorithm based on Lyapunov optimization to calculate them in real time under uncertain workload arrivals and fluctuating grid carbon intensity. Grid operators can issue dispatch signals within the flexibility intervals to fully leverage DC power flexibility. Furthermore, we introduce a fast disaggregation step for the fulfillment of these dispatch signals. Our algorithm is prediction-free, highly adaptable to uncertainties, achieves a bounded worst-case workload request delay, and provides a strict feasibility guarantee for disaggregation. Numerical simulations demonstrate that our approach effectively quantifies and unlocks DC power flexibility while satisfying carbon emission constraints, outperforming baseline algorithms.
\end{abstract}


\IEEEpeerreviewmaketitle

\section{Introduction}
% \IEEEPARstart

Data center (DC) infrastructure is experiencing significant growth, driven by increasing artificial
intelligence (AI) workloads. The global capacity of DCs is projected to double between 2026 and 2030~\cite{jll_data_center_2026}. While DCs offer adequate computational power, their substantial energy consumption not only poses a challenge on the power grid, but also leads to considerable carbon emissions~\cite{bian2024cafe}. Fortunately, DCs can also offer substantial spatio-temporal flexibility to help address these challenges~\cite{zhou2024ai}. Spatial flexibility comes from routing workloads to different DCs for processing, and temporal flexibility comes from deferring delay-tolerant workloads to later periods.

%To alleviate these pressures, researchers are turning to energy management of DCs by strategically scheduling their workloads to unlock the inherent spatio-temporal flexibility of DCs~\cite{zhou2024ai}. %This gives rise to the focus of this paper, which is to manage DC workloads to quantify and then provide demand-side power flexibility to power grids while controlling their carbon emissions.

Extensive research has been devoted to workload scheduling of DCs. 
Reference \cite{CHEN2021116935} explored the potential of DCs to improve the operational flexibility of distribution networks. A planning model considering carbon-oriented demand responses of DC loads was proposed to evaluate the carbon emission reduction potential of DCs \cite{du2025estimating}.
In~\cite{sun2025privacy}, a day-ahead scheduling framework was proposed to coordinate the operation of DCs, which is implemented by a privacy-preserving distributed algorithm. To handle the uncertainty therein, a stochastic optimization model was developed to coordinate the operation of power systems and DCs~\cite{niu2021spacial}.
Moreover, reference \cite{chen2016robust} formulated the workload management problem of DCs as a robust optimization problem to minimize the operational cost.
However, the above methods are all \textit{offline} methods, which depend heavily on future predictive information that is often difficult to obtain accurately. Moreover, these methods often fail to incorporate the sequentially revealed real-time dispatch signals from the power grid operator, and may therefore underestimate or underutilize the power flexibility of DCs.

% However, the aforementioned offline methods suffer from the following limitations:
% \begin{itemize}
%     \item Offline methods depend heavily on future predictive information, which is often difficult to obtain accurately.
%     \item When DCs provide power flexibility to grids, the power grid operator must issue dispatch signals to harness these resources, leading to an interactive process. However, such critical signals are unavailable in an offline setting.
%     \item In an offline environment, workload scheduling decisions cannot be utilized as feedback signals to update DC system dynamics in real time.
% \end{itemize}

To overcome the limitations of offline methods, a body of literature manages DC workloads through online schemes such as model predictive control (MPC).
Reference~\cite{zhang2022exploring} used MPC to manage DC operations where the impact of grid signals such as carbon emission rates was explored.
Reference \cite{Carnerero2022proba} employed MPC to improve the energy efficiency of DCs under temperature and service constraints.
Reference \cite{cupelli2018data} proposed a DC control framework for demand response participation based on MPC. However, the above MPC-based methods still rely on predictive information such as future workloads
and can be computationally demanding. Alternatively, reinforcement learning (RL) is an emerging paradigm for online DC workload scheduling~\cite{KAHIL2025125734}. 
To minimize the energy cost, reference \cite{wei2021deep} utilized RL to jointly control DCs and HVAC systems within buildings.
Reference \cite{ran2023optimizing} proposed a deep RL approach to jointly optimize the workload scheduling and cooling control in a DC.
The work \cite{wei2023multi} developed a RL-based method to achieve balanced resource utilization of DCs. Reference \cite{liu2023online}~designed a deep RL-based  method to handle the large-scale heterogeneous computing tasks in DCs.
Nevertheless, the performance of RL-based schemes heavily rely on sufficient training, and RL can
face difficulties generalizing to unseen scenarios.

Lyapunov optimization (LO) is a prediction-free and training-free approach to deliver online decisions in uncertain environments~\cite{neely2010stochastic}, and thus provides an alternative to prediction-based or learned control policies. LO has been used in various areas including electric vehicle (EV) charging scheduling~\cite{wang2025online,menghwar2024market}, building temperature control~\cite{zheng2014distributed} and energy storage control~\cite{huang2012optimal}.  LO also finds applications in DC workload management problems.
For instance, the work~\cite{gao2024online} designed a hierarchical scheme to coordinate the operation and energy dispatch of DCs, where LO was utilized to perform electricity purchase decisions for DCs. Reference \cite{sun2023battery} designed a distributed online algorithm to minimize the operational cost of DCs based on LO. Moreover, references~\cite{dou2017carbon,lin2023carbon,he2019cost} further address DC workload management with consideration of carbon emissions. Based on LO, reference \cite{dou2017carbon} proposed a carbon-aware online algorithm to schedule deferrable workloads considering both electricity bills and the carbon tax.
Within the LO framework, reference \cite{lin2023carbon} explored the cost and carbon reduction potential of DCs with distributed renewable generation.
Reference \cite{he2019cost} developed an online algorithm based on LO to minimize the cost of DCs under carbon emission budgets.

Although the above studies have demonstrated promising results on online DC workload management based on LO, they focus on optimizing a desired workload flow rather than quantifying the power flexibility of DCs, which requires
characterizing a set of realizable power decisions. %For a flexible resource, a \textit{flexibility interval} is usually employed to describe its power flexibility potential, explaining how its power demand can be adjusted between the upper and lower bounds of the interval upon real-time requests from power grids. Moreover, when quantifying the online power flexibility of DCs, we need to guarantee the \textit{disaggregation feasibility} of the flexibility interval~\cite{wang2025tte,chen2021leveraging}, i.e, any power level within the flexibility interval can be realized by a feasible workload flow within the DC system without violating the operational constraints of DCs.
%This paper aims to address this research gap by quantifying the online power flexibility of DCs under carbon emission limits based on the LO with disaggregation feasibility guarantees.
Existing works on online power flexibility quantification have considered flexible resources such as EVs~\cite{wang2025carbon,yan2024real,Li2024distributionally} and  buildings
\cite{CAI2024105531,GASSER2021116653,de2016quantification}.
However, evaluating the power flexibility of DCs introduces new research challenges: (1) Because DCs possess spatio-temporal flexibility, evaluating their flexibility requires considering not only the impact of current decisions on future time slots, as with EVs, but also the interdependence among geographically distributed DCs. (2) Computing resource management and energy management jointly determine DC power consumption and are tightly coupled through IT load. Flexibility evaluation therefore requires projecting a high-dimensional operational feasible region onto the power space, which is computationally intractable.
    % \item Under time-varying and complex environments, it is difficult to achieve flexibility maximization while ensuring the operational requirements of DC system.
    % \item  Geographically distributed DCs often interact with power grids across multiple electric regions, and the demand for flexibility services in each region's power grid is different and time-varying. Consequently, systematically coordinating simultaneous interactions across these distinct electric regions to provide power flexibility presents a formidable challenge.
  %  \item  Carbon intensity exhibits complex dynamics across both spatial and temporal dimensions. Hence, effectively exploring the power flexibility potential of DCs under carbon emission constraints remains non-trivial, as it requires the trade-offs between flexibility provision and workload management.
    %\item Online environments require fast responses to dispatch signals.  Given  the complex operational constraints and spatio-temporal characteristics of DC systems, designing fast algorithms of flexibility quantification and disaggregation suitable for online implementation with a disaggregation feasibility guarantee, is a highly non-trivial task.   
    % For the disaggregation process, we must guarantee feasibility without violating complex constraints of DCs, while simultaneously ensuring fast response to dispatch signals in an online environment. Thus, how to devise a fast  disaggregation method suitable for online implementation with a disaggregation feasibility guarantee, is a non-trivial task. 

We aim to address the above challenges by our proposed online carbon-aware DC power flexibility quantification framework. Our main contributions are two-fold:

%\item  We propose an effective online framework for carbon-aware DC power flexibility quantification. In our framework, we quantify DC power flexibility by power flexibility intervals, which consists of dynamical upper and lower bounds of power. In our proposed framework, we design separate workload flows within DCs for deriving the upper and lower bounds, and such a design also ensures the disaggregation feasibility of flexibility intervals. Then delay virtual queues are constructed to constrain the lower bound, thereby limiting the worst-case workload request delay, while a carbon virtual queue is built to constrain the upper bound for controlling the carbon emission rate of the entire DC system. With these queues, we systematically achieve the DC power flexibility quantification by formulating a stochastic programming problem, which aims to maximize the time-average value of DC power flexibility subject to queue stability conditions.

(1) \textit{\textbf{DC Power Flexibility Quantification Model}}: We propose an optimization model to characterize DC power flexibility through a time-decoupled representation: the Cartesian product of adjustable power demand intervals across time slots. Unlike the high-dimensional polyhedral representations in previous work, the proposed formulation provides a simple, operator-friendly way for incorporating DC flexibility into power system operation while preserving intertemporal operational feasibility. According to the main factors that constrain DC operation, the lower bound of each interval is set by constraints on workload flows and processing requirements while the upper bound is set by constraints on cumulative carbon emissions and their limits. Proposition \ref{Prop-1} proves that the proposed model guarantees disaggregation feasibility, i.e., any power trajectory the system operator picks within these intervals can be achieved by feasible DC operation strategies.

(2) \textit{\textbf{Online Algorithm}}: We develop an online algorithm to solve the proposed DC power flexibility quantification model using LO techniques. A key innovation is the joint design of virtual queues that capture intertemporal coupling in both DC workload flows and carbon emissions, two factors rarely addressed in existing online algorithms. Moreover, we develop an integrated disaggregation and feedback mechanism that translates the power system operator's dispatch commands into real-time DC operating decisions and feeds the resulting operational states back into subsequent flexibility quantification. This closed-loop design explicitly accounts for how the realization of current flexibility shapes future flexibility, linking online flexibility quantification with real-time operational implementation.

The rest of this paper is organized as follows. Section \ref{sec:formulation} introduces the DC power flexibility quantification model, which is solved by an online algorithm in Section \ref{sec:online_algorithm}.
The numerical simulations are presented in Section \ref{sec:case} and Section \ref{sec:conclusion} concludes this work.

\section{Problem Formulation}
\label{sec:formulation}
In this section, we first introduce the system framework. Then we model the DC system and formulate an offline problem to quantify the DC power flexibility. 

\subsection{System framework}
\begin{figure}
    \centering
    \includegraphics[width = 0.9\linewidth]{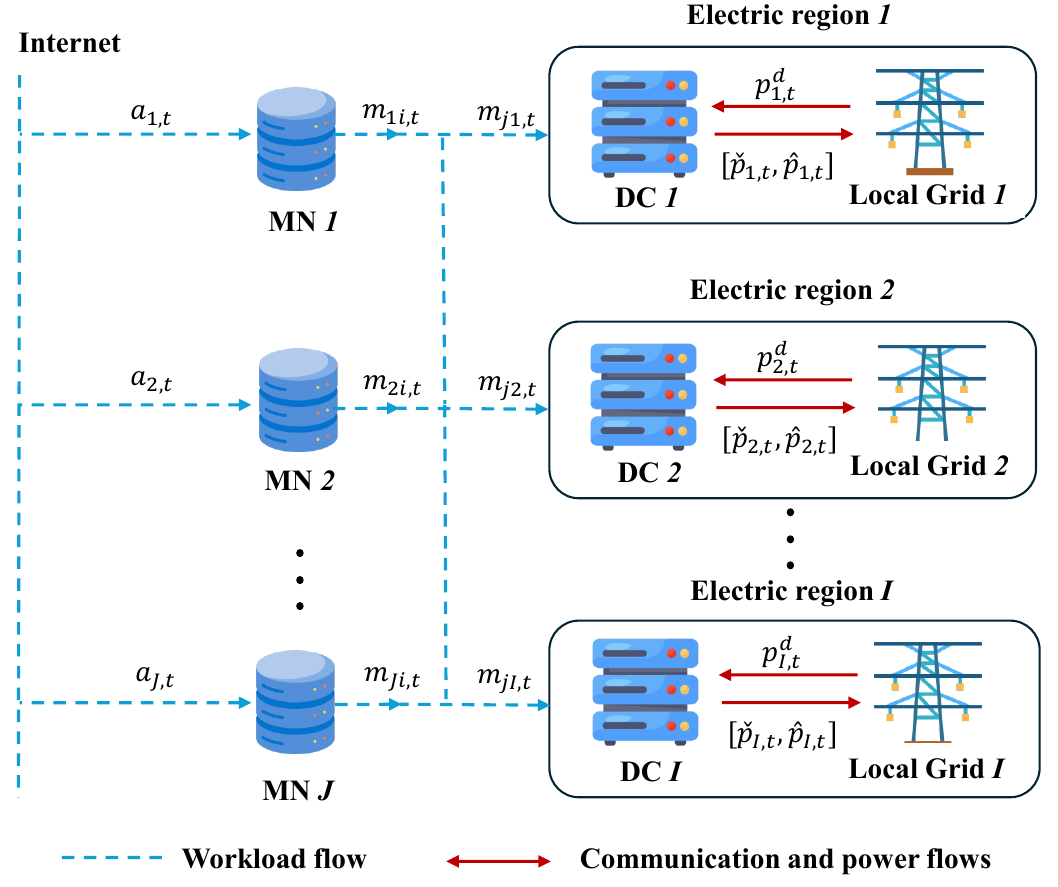} 
    % \vspace{-1em}
    \caption{ Geographically distributed DC system.}
    \label{fig:system}
\end{figure}

We consider a DC system in Figure~\ref{fig:system} with a set of mapping nodes (MNs) $\mathcal{J} := \{1,2,...,J\} $ and a set of geographically distributed DCs $\mathcal{I} := \{1,2,...,I\}$. The time slot is indexed by $t$ and the entire system operates in a discrete time horizon $\mathcal{T} := \{1,2,...,T\}$. The MNs, which serve as the front-end components, collect workload requests from users, and forward the requests to back-end DCs for processing. We use $\Omega_j^f$  to denote the set of back-end DCs connected to the front-end MN $j \in \mathcal{J}$, and use $\Omega_i^b$ to denote the set of front-end MNs connected to the back-end DC $i \in \mathcal{I}$. We assume that the electricity consumption of a unit of workload is a constant so that we can quantify the workload by the amount of electricity consumption. Workloads are processed in a first-in-first-out (FIFO) order: workload arriving earlier at a DC or MN
has priority for processing or forwarding, respectively.  The operation of the entire DC system is jointly accomplished by the MNs and DCs, where the front-end MNs are responsible for workload distribution and the back-end DCs process workload requests by purchasing electricity from local power grids.

In this work, we focus on delay-tolerant workloads, and our proposed framework can also incorporate delay sensitive workloads by assigning them higher processing priority.  The MNs have the flexibility to choose when to transfer workloads to geographically distributed DCs. The DCs have the flexibility to determine when to process the workloads received from MNs. Therefore, the power consumption of each DC $i$ has certain flexibility formally defined below.
%At each time slot $t$, the power level of each distributed DC $i$ can be adjusted within a DC power flexibility interval $[\check{p}_{i,t},\hat{p}_{i,t}]$, which is defined below.
% \begin{definition}(DC power flexibility interval).\label{def-1}
%     $[\check{p}_{i,t},\hat{p}_{i,t}]$ is called the \textbf{power flexibility interval} of DC $i$ at time $t$ if for any $p_{i,t}^d \in [\check{p}_{i,t},\hat{p}_{i,t}]$, there exists a feasible workload dispatch and processing strategy for the DC system such that the power consumption of DC $i$ is $p_{i,t}^d$.
% \end{definition}
\begin{definition}(DC power flexibility intervals)\footnote{This definition concerns one-slot feasibility, not
the feasibility of every future sequence of interval selections.}: At time $t$, let $\mathcal{F}_t := \Pi_{i \in \mathcal{I}} [\check{p}_{i,t},\hat{p}_{i,t}]$ represent the \textit{\textbf{DC power flexibility intervals}}, which satisfy that for every $p_t^d \in \mathcal{F}_t$, there exists a simultaneous workload routing
decision satisfying the current-slot link, workload-availability, and processing-capacity constraints, with DC $i$ consuming $p_{i,t}^d$ for every $i$.  \label{def-1}
\end{definition}
% For each DC, the DC operator quantifies such power flexibility interval, and reports it to the corresponding local power grid for flexibility service provision. After receiving the flexibility interval, the local grid operator can choose any power level within the flexibility interval by sending a dispatch signal (i.e., $p_{i,t}^d \in [\check{p}_{i,t},\hat{p}_{i,t}]$ for DC $i$) to the DC operator to utilize this flexibility asset.
Intuitively, the lower bound $\check{p}_{i,t}$ of the power flexibility interval is constrained by workload processing requirements, and thus when defer processing tasks, the DC should guarantee their timely completion. In contrast, beyond the constraint of workload processing capacity, the upper bound $\hat{p}_{i,t}$ is also constrained by carbon emission limits, capping DC power consumption to ensure compliance with emission targets. As the carbon emissions currently represent a greater concern for DCs, here the carbon emission constraint is assumed to be tighter than the processing capacity constraint.

% {\color{red}(What if the upper bound determined by the maximum IT workload processing power is less than that determined by carbon limits?)}

%related to deferrable nature of the workload, allowing the DC operator to decrease the power level by delaying the load's completion. Conversely, the upper bound $\hat{p}_{i,t}$ represents the goal of increasing the power level as much as possible, but it is constrained by the decarbonization goal since the demand-side DCs consume a large amount of energy, resulting in massive carbon emissions.

In the following, we will introduce the model to quantify the DCs' power flexibility considering carbon emission limits. We use $ \check{(\cdot)}$ and $ \hat{(\cdot)}$ to denote the variables associated with the lower and upper bounds of the power flexibility interval, respectively.

%model the entire DC system and formulate the flexibility quantification problem of DCs under carbon emission limits.

\subsection{Data center flexibility quantification}
 %Two workload flows are designed for the upper and lower bounds.
\subsubsection{Constraints for the lower bound of power flexibility intervals}
We begin by defining the constraints related to the lower bound $\check{p}_{i,t}$. Let $a_{j,t}$ denote the amount of the workload that arrives at the front-end MN $j$ at time $t$. After receiving the workload, the MN $j$ allocates the workload to the back-end DC $i \in \Omega_j^f$ for processing. Let $\check{m}_{ji,t}$ be the workload transferred from MN $j$ to DC $i$ related to the lower bound. The remaining workload is buffered in a corresponding queue for each MN and the queue at MN $j$ evolves with dynamics:
\begin{equation}\label{eq:queue_f}
Q_{f,j,t+1} = Q_{f,j,t} + a_{j,t} - \sum\nolimits_{i \in \Omega_{j}^f} \check{m}_{ji,t},
\end{equation}
where $Q_{f,j,t}$ is the queue length of MN $j$ at time $t$.

The back-end DC $i$ receives workload from MN $j \in \Omega_i^b$ and draws energy from the power grid. Each back-end DC also holds a queue to store the unserved workloads. We use $Q_{b,i,t}$ to denote the queue length of DC $i$ at time $t$, whose dynamics is: 
\begin{equation}\label{eq:queue_b}
Q_{b,i,t+1} = Q_{b,i,t} + \sum\nolimits_{j \in \Omega_{i}^b} \check{m}_{ji,t} - \check{p}_{i,t}.
\end{equation}
Here, $\check{p}_{i,t}$ represents the amount of energy DC $i$ bought from the power grid at time $t$ for workload processing. For simplicity, we assume the time interval is 1 hour, so energy and power can be used equivalently.

The following constraint ensures that the amount of the workload transferred between MN $j$ and DC $i$ is upper bounded by the link bandwidth $M_{ji}$, measured in
MWh-equivalent workloads per slot:
\begin{equation}\label{eq:lower_link_bandwidth}
    0 \leq \check{m}_{ji,t} \leq M_{ji}, \forall j, \forall i \in \Omega_j^f, \forall t.
\end{equation}

The lower bound of the power flexibility interval at time $t$ is limited by
\begin{equation}
    0 \leq \check{p}_{i,t} \leq \check{p}_{i,t}^{max}, \forall i,t,
\end{equation}
where $\check{p}_{i,t}^{max}$ is determined by both the workload processing capacity $\bar{P}_i$ and the available workload in the queue $Q_{b,i,t}$.
\begin{equation}\label{eq:lower_bound_limit}
    \check{p}_{i,t}^{max} = \min \{\bar{P}_{i}, Q_{b,i,t} + \sum\nolimits_{j \in \Omega_{i}^b} \check{m}_{ji,t}\}, \forall i,t.
\end{equation}
% {\color{red}(Given this constraint (5), it seems that we do not need $\max\{xx,0\}$ in (2))}

%Unrestricted lowering of the bound $\check{p}_{i,t},\forall t$ would potentially result in  infinite unserved workloads in the queue. Therefore, stability in queue length should be carefully maintained over the time horizon, i.e.,
To avoid unlimited unserved workloads, we seek the stability of workload queues:
\begin{equation}\label{eq:workload_bound}
    \lim_{T \to \infty}\frac{1}{T}\sum_{t \in \mathcal{T}} \{ \sum_{j \in \mathcal{J}}\mathbb{E}[Q_{f,j,t}] + \sum_{i \in \mathcal{I}}\mathbb{E}[Q_{b,i,t}] \}  < \infty.
\end{equation}

\subsubsection{Constraints for the upper bound of power flexibility intervals} We then construct the constraints for the upper bound $\hat{p}_{i,t}$ by adopting a separate set of variables to describe the flow of workloads. The transferable workload in the flow of the upper bound, denoted by $\hat{m}_{ji,t}$, satisfies the following:  
\begin{equation}\label{eq:upper_link_bandwidth}
     \check{m}_{ji,t} \leq \hat{m}_{ji,t} \leq M_{ji}, \forall j,i,t.
\end{equation}
%This means that $\hat{m}_{ji,t}$ is also physically constrained by $M_{ji,t}$. 
For MN $j$ at time $t$, the sum of the workloads allocated to different DCs cannot exceed the available workload:
\begin{equation}\label{eq8}
    \sum\nolimits_{i \in \Omega_{j}^f} \hat{m}_{ji,t} \leq Q_{f,j,t} + a_{j,t}.
\end{equation}
Equations \eqref{eq:upper_link_bandwidth}-\eqref{eq8} also enforce $\sum\nolimits_{i \in \Omega_{j}^f} \check{m}_{ji,t} \le  Q_{f,j,t}+a_{j,t}$, so the front-end queue backlog in \eqref{eq:queue_f} is nonnegative.

Similarly, we impose the following constraints on $\hat{p}_{i,t}$:
\begin{equation}
    \check{p}_{i,t} \leq \hat{p}_{i,t} \leq \hat{p}_{i,t}^{max}, \forall i,t.
\end{equation}

To align with the workload flow of the upper bound, $\hat{m}_{ji,t}$ is used to determine $\hat{p}_{i,t}^{max}$.
\begin{equation}\label{eq:upper_bound_limit}
    \hat{p}_{i,t}^{max} = \min \{\bar{P}_{i}, Q_{b,i,t} + \sum\nolimits_{j \in \Omega_{i}^b} \hat{m}_{ji,t}\}, \forall i,t.
\end{equation}

The long-term time-average carbon emission limit of the entire DC system is described by
\begin{equation}\label{eq:carbon_cons}
    \lim_{T \to \infty}\frac{1}{T}\sum\nolimits_{t \in \mathcal{T}} \sum\nolimits_{ i \in \mathcal{I}}\mathbb{E} [\beta_{i,t}\hat{p}_{i,t}] \leq E,
\end{equation}
where $\beta_{i,t}$ is the time-varying carbon intensity of the electricity that DC $i$ draws from the power grid at time $t$. Note that the value $\beta_{i,t}$ varies across different DCs, since the DCs are located in various electric regions, each with a distinct energy mix. $E$ is the long-term time-average carbon allowance, which can be estimated by averaging the total carbon emission requirement throughout the time horizon.

\subsubsection{DC power flexibility quantification model}
With these constraints, we formulate the following model to maximize the time-average total power flexibility of the DC system:
\begin{subequations}
    \begin{align}
 \mathbf{P1\colon} \; \max\; &  \lim_{T \to \infty}\frac{1}{T} \sum_{t =1}^T   \mathbb{E}\left[\sum_{i \in \mathcal{I}} \gamma_{i,t}(\Hat{p}_{i,t} -\Check{p}_{i,t})\right], \label{eq:P1_a}\\
&\textrm{s.t.} \qquad (\text{1}) - (\text{11}),
    \end{align}
\end{subequations}
where $\gamma_{i,t}$ serves as a reward signal issued by the power grid operator for back-end DC $i$'s flexibility provision. %$\gamma_{i,t}$ can be time-varying and is regulated in accordance with local power grid conditions to guide the flexibility quantification. 
If $\gamma_{i,t}$ is set to an identical constant for all back-end DCs, then $\mathbf{P1}$ aims to maximize
the average aggregate interval width, without considering the grid-side preferences.
% {\color{red}(If all $\gamma_{i,t}$ are equal, instead of $\gamma_{i,t}=1$; For example, $\gamma_{i,t}=0.5,\forall i,t$ is equivalent to $\gamma_{i,t}=1,\forall i,t$)}

%It is not a suitable choice to solve $\mathbf{P1}$ directly, and this can be explained by two reasons.
However, it is difficult to solve $\mathbf{P1}$ directly since it requires accurate information about the future stochastic environment, which is rarely available in practice. Moreover, this static model tends to underestimate DC power flexibility because grid dispatch signals unfold sequentially and could influence the future DC power flexibility. Therefore, an online algorithm for solving $\mathbf{P1}$ is needed.

%the future workloads and grid carbon intensity across various locations are uncertain and hard to predict accurately. {Second, as $\mathbf{P1}$ is an offline formulation,  when quantifying the flexibility, the impact of real-time dispatch signal $p_{i,t}^d$ on DC flexibility characterization is not considered.}
%The above reasons call for an efficient online algorithm to tackle $\mathbf{P1}$ since in an online mode, uncertain factors and dispatch signals would be sequentially revealed.

%In the next section, we will develop an online algorithm based on LO to address the above issues.

\section{Online Algorithm}
\label{sec:online_algorithm}
%In this section, we design a set of virtual queues, and reformulate the DC power flexibility quantification problem with these virtual queues. Supported by LO, we then develop an online algorithm to solve the problem. Finally, we further design a disaggregation step and update the queues with dispatch signals issued by grid operators.
In this section, we develop an online algorithm for solving the DC power flexibility quantification model $\mathbf{P1}$ based on LO, and design disaggregation and update steps to account for the impact of the latest dispatch signal on future flexibility.

%There are two issues that make it challenging to design an effective online algorithm. 1) Although workloads are delay tolerant, the imposition of rigorous bounds on maximum workload request delay is a non-trivial task. 2) The decision variables (e.g., $\hat{p}_{i,t}$) in the problem are time-coupled due to the existence of the long-term carbon compliance constraint \eqref{eq:carbon_cons}. To address these issues, we introduce a set of virtual queues in the following development.

\subsection{Design of virtual queues}
First, to bound the worst-case workload processing delay, we design virtual delay queues $Z_{f,j,t}$ and $Z_{b,i,t}$ for workload queues $Q_{f,j,t}$ and $Q_{b,i,t}$, respectively. These virtual delay queues evolve with dynamics as follows:
\begin{equation}\label{eq:queue_Z_f}
Z_{f,j,t+1} = \max\{Z_{f,j,t} - \sum\nolimits_{i \in \Omega_{j}^f}\check{m}_{ji,t} + \delta_{f,j} \mathbbm{1}_{Q_{f,j,t}},0 \},
\end{equation}
\begin{equation}\label{eq:queue_Z_b}
Z_{b,i,t+1} = \max\{Z_{b,i,t} - \check{p}_{i,t} + \delta_{b,i} \mathbbm{1}_{Q_{b,i,t}},0 \},
\end{equation}
where $\delta_{f,j}$ and $ \delta_{b,i}$ are predefined constants to control the growth rate of the virtual queues. $\mathbbm{1}_y$ is an indicator function with $\mathbbm{1}_y = 1$ if $y > 0$ and $\mathbbm{1}_y = 0$ otherwise. We set initial lengths of the virtual delay queues to be 0.  %One can observe that the above virtual delay queues interact with the workload flow of the lower bound. This is because, as mentioned previously, the lower bound of the power flexibility interval is mainly constrained by the workload requirements. 
In the following theorem, we will show that the integration of these virtual delay queues provides an upper bound of the worst-case workload processing delay.

\begin{theorem}\label{theorem_1}
Suppose requests are served in FIFO order, $\delta_{f,j}>0$ and $\delta_{b,i}>0$, and the actual queues updated by \eqref{eq:update_1}-\eqref{eq:update_4} obey, for all slots,
$Q_{f,j,t} \leq Q_{f,j}^{max}$, $Q_{b,i,t} \leq Q_{b,i}^{max}$, $Z_{f,j,t} \leq Z_{f,j}^{max}$ and $Z_{b,i,t} \leq Z_{b,i}^{max}$,
with finite deterministic constants. A conservative end-to-end
delay bound for work entering MN $j$ is $D_{f,j} + \max_{i \in \Omega_j^f}D_{b,i}$, where the bounds in slots are $D_{f,j} = 1+\lceil Q_{f,j}^{max} + Z_{f,j}^{max}\rceil /{\delta_{f,j}}$ and $D_{b,i} = 1+ \lceil Q_{b,i}^{max} + Z_{b,i}^{max}\rceil /{\delta_{b,i}}$.
\end{theorem}

%\begin{theorem}\label{theorem_1}
%Suppose that workload requests are served in FIFO order and that, for all $t \in \mathcal{T}$, $Q_{f,j,t} \leq Q_{f,j}^{max}$, $Q_{b,i,t} \leq Q_{b,i}^{max}$, $Z_{f,j,t} \leq Z_{f,j}^{max}$ and $Z_{b,i,t} \leq Z_{b,i}^{max}$. 
  %  Assume $Q_{f,j,t}$, $Q_{b,i,t}$, $Z_{f,j,t}$ and $Z_{b,i,t}$ can be ensured to be upper bounded for all time slots $t \in \mathcal{T}$ by an algorithm, i.e., $Q_{f,j,t} \leq Q_{f,j}^{max}$, $Q_{b,i,t} \leq Q_{b,i}^{max}$, $Z_{f,j,t} \leq Z_{f,j}^{max}$ and $Z_{b,i,t} \leq Z_{b,i}^{max}$. Suppose workload requests are processed in FIFO manner. 
 % The worst-case queuing delay of the request that arrives at front-end MN $j$ is bounded by $D_{f,j} + \max_{i \in \Omega_j^f}D_{b,i}$, where $D_{f,j} = (Q_{f,j}^{max} + Z_{f,j}^{max}) /{\delta_{f,j}}$ and $D_{b,i} = (Q_{b,i}^{max} + Z_{b,i}^{max}) /{\delta_{b,i}}$.
%\end{theorem}
The proof of Theorem \ref{theorem_1} can be found in Appendix \ref{app_A}.
Furthermore, we construct another  virtual queue for  the time-coupled carbon-related constraint \eqref{eq:carbon_cons}:
\begin{equation}\label{eq:carbon_queue}
    Q_{c,t+1} = \max\{Q_{c,t} + \sum\nolimits_{i \in \mathcal{I}} \beta_{i,t} \hat{p}_{i,t} - E, 0\},
\end{equation}
where $ Q_{c,t}$ is the backlog of the virtual carbon queue, recording the time-accumulated difference between $ \sum\nolimits_{i \in \mathcal{I}} \beta_{i,t} \hat{p}_{i,t}$ and $E$. Its mean-rate stability is a sufficient condition
for \eqref{eq:carbon_cons}:
%With this virtual carbon queue, the constraint \eqref{eq:carbon_cons} can be relaxed into the mean-rate-stable constraint of the virtual queue, i.e., 
\begin{equation}\label{eq:carbon_queue_stab}
    \lim_{T \to \infty} \frac{\mathbb{E}[Q_{c,T+1}]}{T} = 0.
\end{equation}
This is because, by \eqref{eq:carbon_queue}, we have:
%With this virtual carbon queue, we are able to decompose the time-coupled carbon constraint  into each time slot, since we only need to focus on maintaining the stability of the carbon virtual queue at each time slot to meet the constraint \eqref{eq:carbon_cons}. In the following, we show how this stability ensures the satisfaction of the related carbon emission constraint. By \eqref{eq:carbon_queue}, we have:
\begin{equation}
    Q_{c,t+1} \geq Q_{c,t} + \sum\nolimits_{i \in \mathcal{I}} \beta_{i,t}\hat{p}_{i,t} - E.
\end{equation}
Summing over slots, taking expectations, and dividing by $T$
gives the following
\begin{equation}
    \frac{\sum_{t=1}^{T} \sum_{i \in \mathcal{I}}\mathbb{E}[\beta_{i,t}\hat{p}_{i,t} - E]}{T} \le \frac{\mathbb{E}[Q_{c,T+1}] - \mathbb{E}[Q_{c,1}] }{T} =0,
\end{equation}

% The above inequality indicates that constraint \eqref{eq:carbon_cons} can be satisfied under the following queue stability condition:

% \begin{equation}\label{eq:carbon_queue_stab}
%     \lim_{T \to \infty} \frac{\mathbb{E}[Q_{c,T+1}]}{T} = 0.
% \end{equation}

 With the virtual queues above, we reformulate the problem $\mathbf{P1}$  as follows:
\begin{subequations}
    \begin{align}
 \mathbf{P2\colon} & \;\min\;   \lim_{T \to \infty}\frac{1}{T} \sum_{t =1}^T  \mathbb{E}[ -\sum_{i \in \mathcal{I}} \gamma_{i,t}(\Hat{p}_{i,t} -\Check{p}_{i,t})], \label{eq:P2_a}\\
&\textrm{s.t.} \quad \eqref{eq:queue_f} - \eqref{eq:upper_bound_limit}, \eqref{eq:queue_Z_f}-\eqref{eq:carbon_queue}, \\
&  Z_{f,j,t},  Z_{b,i,t}, Q_{c,t}\text{ are mean rate stable}.
    \end{align}
\end{subequations}
%Note that in the objective function, we utilize the fact that maximizing the average power flexibility is equivalent to minimizing its negative form. 
However, $\mathbf{P2}$ is still an offline optimization model. To address this issue, we apply the LO technique to further turn it into an online counterpart.
%However, we have so far constructed a set of virtual queues to address the underlying issues mentioned earlier. In the next subsection, we will enter into the landscape of LO and explore how these virtual queues help to develop an online algorithm to solve $\mathbf{P2}$ with theoretical guarantees.

\subsection{Lyapunov optimization}
Let $\mathbf{{\Theta}}_t = [{Q}_{f,j,t},{Z}_{f,j,t}, \forall j;{Q}_{b,i,t},{Z}_{b,i,t},\forall i;Q_{c,t}]$ be  the vector collecting the information of all queues at time $t$.
The Lyapunov function is defined to measure the aggregated queue backlog of the entire DC system:
\begin{equation}\label{eq:lay_func}
    L(\mathbf{{\Theta}}_t) \stackrel{\triangle}{=} \frac{1}{2}\alpha Q_{c,t}^2 + \frac{1}{2}\sum_{j \in \mathcal{J}}[Q_{f,j,t}^2 + Z_{f,j,t}^2] + \frac{1}{2}\sum_{i \in \mathcal{I}}[Q_{b,i,t}^2 + Z_{b,i,t}^2],
\end{equation}
where $\alpha$ is a non-negative parameter used to regulate the importance of carbon considerations in the Lyapunov function.

The Lyapunov drift is given by the expectation of the one-slot difference in the Lyapunov function:
{\begin{equation}
    \Delta (\mathbf{{\Theta}}_t) \stackrel{\triangle}{=} \mathbb{E}[L(\mathbf{{\Theta}}_{t+1}) - L(\mathbf{{\Theta}}_t) \mid \mathbf{{\Theta}}_t ].
\end{equation}}

A large value of $L(\mathbf{{\Theta}}_t)$ reflects that queues are experiencing  congestion, potentially caused by either substantial unserved workloads or significant carbon emissions. Minimizing $ \Delta (\mathbf{{\Theta}}_t)$ pushes the queue states towards less congested, maintaining the stability of all queues and thereby ensuring the satisfaction of the related long-term queue stability constraints. However, we can not greedily minimize $\Delta (\mathbf{{\Theta}}_t)$ at each time slot $t$, as this may contribute to a low power flexibility level. Therefore, we build the drift-plus-penalty term to balance the Lyapunov drift and the expected power flexibility at time $t$:
\begin{equation}\label{eq:drift_plus_penalty}
  \Delta (\mathbf{{\Theta}}_t) + V \mathbb{E}[- \sum\nolimits_{i \in \mathcal{I}} \gamma_{i,t}(\Hat{p}_{i,t} -\Check{p}_{i,t}) ~| ~{\Theta}_t],
\end{equation}
where $V$ is a non-negative parameter, balancing the queue stability constraints and the objective function.
% Directly minimizing the drift-plus-penalty term in \eqref{eq:drift_plus_penalty} usually requires future information
%Due to the time-coupling setting in the Lyapunov drift $\Delta (\mathbf{{\Theta}}_t)$, 
To further remove the time-coupling terms in the Lyapunov drift, 
we actually build an upper bound of $\Delta (\mathbf{{\Theta}}_t)$ to design our online algorithm as follows.

% {\color{red}(Since you have changed constraint (1), please update the following derivations.)}

From \eqref{eq:queue_f}, we can have the following:
\begin{equation} \label{eq:inequality_bound}
    Q_{f,i,t+1}^2 \leq (Q_{f,i,t} + a_{j,t} - \sum\nolimits_{i \in \Omega_j^f}\check{m}_{ji,t})^2.
\end{equation}

Then based on \eqref{eq:inequality_bound}, we can have the following inequality:
 \begin{align*}
        &\frac{1}{2} \sum_{j \in \mathcal{J}} (Q_{f,j,t+1}^2 - Q_{f,j,t}^2) \leq \frac{1}{2}\sum_{j \in \mathcal{J}}(a_{j,t} - \sum\nolimits_{i \in \Omega_j^f}\check{m}_{ji,t})^2 \\
        & + Q_{f,j,t}(a_{j,t} - \sum\nolimits_{i \in \Omega_j^f}\check{m}_{ji,t}).\stepcounter{equation}\tag{\theequation}
\end{align*}

From \eqref{eq:queue_b}, we have:
\begin{align*}
        &\frac{1}{2} \sum_{i \in \mathcal{I}} (Q_{b,i,t+1}^2 - Q_{b,i,t}^2) \leq \frac{1}{2}\sum_{i \in \mathcal{I}}(\sum\nolimits_{j \in \Omega_i^b}\check{m}_{ji,t} - \check{p}_{i,t})^2\\
        & + Q_{b,i,t}(\sum\nolimits_{j \in \Omega_i^b}\check{m}_{ji,t} - \check{p}_{i,t}).\stepcounter{equation}\tag{\theequation}
\end{align*}

Equation \eqref{eq:queue_Z_f} and $\mathbb{I}_{Q_{f,j,t}} \le 1$ yield
%implies the fact that $Z_{f,j,t+1} \leq \max[Z_{f,j,t} - \sum_{i \in \Omega_{j}^f}\check{m}_{ji,t} + \delta_{f,j},0]$, and thus for the virtual delay queue $ Z_{f,j,t}$, we have:
\begin{align*}
        &\frac{1}{2} \sum_{j \in \mathcal{J}} (Z_{f,j,t+1}^2 - Z_{f,j,t}^2) \leq \frac{1}{2}\sum_{j \in \mathcal{J}} ( \delta_{f,j}- \sum\nolimits_{i \in \Omega_j^f}\check{m}_{ji,t})^2\\
        & + Z_{f,i,t}(\delta_{f,j} - \sum\nolimits_{i \in \Omega_j^f}\check{m}_{ji,t}).\stepcounter{equation}\tag{\theequation}
\end{align*}

In a similar way, we have
\begin{equation}
    \begin{aligned}
        &\frac{1}{2} \sum_{i \in \mathcal{I}} (Z_{b,i,t+1}^2 - Z_{b,i,t}^2) \leq \frac{1}{2}\sum_{i \in \mathcal{I}}( \delta_{b,i}- \check{p}_{i,t})^2\\
        & + Z_{b,i,t}(\delta_{b,i} - \check{p}_{i,t}).
    \end{aligned}
\end{equation}

For the virtual carbon queue, we have:
\begin{align*}\label{Q_c_bound}
    &\quad \;\frac{1}{2}(Q_{c,t+1}^2 - Q_{c,t}^2)\stepcounter{equation}\tag{\theequation}\\
    & \leq \frac{1}{2} ( \sum\nolimits_{i \in \mathcal{I}} \beta_{i,t}\hat{p}_{i,t} - E)^2 +  Q_{c,t}(\sum\nolimits_{i \in \mathcal{I}} \beta_{i,t}\hat{p}_{i,t} - E).
\end{align*}

Summing these bounds gives
\begin{align*}\label{eq:upper_bound}
         &\Delta (\mathbf{{\Theta}}_t) -  V \mathbb{E}[\sum_{i \in \mathcal{I}} \gamma_{i,t}(\hat{p}_{i,t}- \check{p}_{i,t})] \leq   - V  \mathbb{E}[ \sum_{i \in \mathcal{I}}\gamma_{i,t}(\hat{p}_{i,t}-\check{p}_{i,t})]\\
         & +\frac{1}{2}\sum_{j \in \mathcal{J}}\bigg[(a_{j,t} - \sum\nolimits_{i \in \Omega_j^f}\check{m}_{ji,t})^2 +  ( \delta_{f,j}- \sum\nolimits_{i \in \Omega_j^f}\check{m}_{ji,t})^2 \bigg] \\
         & +\frac{1}{2}\sum_{i \in \mathcal{I}}\bigg[(\sum\nolimits_{j \in \Omega_i^b}\check{m}_{ji,t} - \check{p}_{i,t})^2 +  ( \delta_{b,i}- \check{p}_{i,t})^2 \bigg]\stepcounter{equation}  \\
         & +\sum_{j \in \mathcal{J}}\bigg[Q_{f,j,t}a_{j,t} + Z_{f,j,t}\delta_{f,j} - (Q_{f,j,t} + Z_{f,j,t})\sum_{i \in \Omega_j^f}\check{m}_{ji,t} \bigg]\\
         & + \sum_{i \in \mathcal{I}}\bigg[Q_{b,i,t}\sum_{j \in \Omega_i^b}\check{m}_{ji,t} + Z_{b,i,t}\delta_{b,i} - (Q_{b,i,t} + Z_{b,i,t})\check{p}_{i,t} \bigg]\\
         & + \frac{1}{2} \alpha ( \sum\nolimits_{i \in \mathcal{I}} \beta_{i,t}\hat{p}_{i,t} - E)^2 + \alpha Q_{c,t}(\sum\nolimits_{i \in \mathcal{I}} \beta_{i,t}\hat{p}_{i,t} - E). \tag{\theequation}
\end{align*}

\subsection{Real-time algorithm}
By minimizing the right-hand side of \eqref{eq:upper_bound}, we have the following online optimization problem:
\begin{subequations}
    \begin{align}
 &\mathbf{P3\colon} \quad \min \quad  g_t, \label{eq:P3_a}\\
\textrm{s.t.} &\quad \eqref{eq:lower_link_bandwidth} - \eqref{eq:lower_bound_limit}, \eqref{eq:upper_link_bandwidth} - \eqref{eq:upper_bound_limit},
    \end{align}
\end{subequations}
where $g_t$ is obtained by dropping the constant terms in the right-hand side of \eqref{eq:upper_bound}. Note that $\mathbf{P3}$ is a convex optimization, which can be easily tackled by optimization solvers. %At each time $t$, the DC operator observes the current queue states $\mathbf{{\Theta}}_t$ and other system states including the arrival rate of the workload of each MN $a_{j,t}$, the grid carbon intensity of each DC location $\beta_i$, and the reward signal $\gamma_{i,t}$ issued by the power grid operator to quantify the real-time power flexibility $[\check{p}_{i,t},\hat{p}_{i,t}]$ across various regions by solving $\mathbf{P3}$. Then the flexibility intervals are reported to the power grid operator who will determine the dispatch setpoints within them.
Moreover, $\mathbf{P3}$ does not rely on future information, thus it can be solved in real time.
%Such real-time decisions are made without future information, allowing the algorithm to be prediction-free.
Furthermore, we provide a theoretical guarantee for the performance of our algorithm. %Our goal is to quantify and then provide demand-side flexibility to the power grid, and thus we incorporate this goal into the the drift-plus-penalty term \eqref{eq:drift_plus_penalty}, where the parameter $V$ controls the impact of power flexibility on the objective function of $\mathbf{P3}$. The following theorem will reveal this intuition by bounding the time-average power flexibility achieved by $\mathbf{P3}$ within a gap from the theoretical optimum under the influence of $V$.

\begin{theorem}\label{theorem_2}
   Let $f^\ast$ and $f^\ast_1$ be the time-average expectation of DC power flexibility obtained by solving $\mathbf{P3}$ and $\mathbf{P2}$, respectively. Then we have $f_1^\ast - f^\ast \leq \frac{B}{V}$, 
   % \begin{equation}
   %     f_1^\ast - f^\ast \leq \frac{B}{V}.
   % \end{equation}
   where $B$ is defined as
\begin{align*}
        B &=  \frac{1}{2}\sum\nolimits_{j \in \mathcal{J}}\big[\max\{A_{j}^2 , \sum_{i \in \Omega_j^f}M_{ji}^2\} + \max \{ \delta_{f,j}^2, \sum_{i \in \Omega_j^f}M_{ji}^2\} \big] \\
         & + \frac{1}{2}\sum\nolimits_{i \in \mathcal{I}}\big[\max\{\sum_{j \in \Omega_i^b}M_{ji}^2 , \check{P}_{i}^2\} +  \max\{ \delta_{b,i}^2, \check{P}_{i}^2\} \big]\\
         & + \sum\nolimits_{j \in \mathcal{J}}\big[ Z_{f,j}^{max}\delta_{f,j} \big] + \sum_{i \in \mathcal{I}}\big[Z_{b,i}^{max}\delta_{b,i} \big] \stepcounter{equation}\tag{\theequation}\\
         & +\frac{1}{2} \alpha \max\{ \sum\nolimits_{i \in \mathcal{I}} \bar{\beta}_{i}^2\hat{P}_{i}^2, E^2\},
\end{align*}
and $A_j$, $\check{P}_i$, $\hat{P}_i$ and $\bar{\beta}_i$ are the upper bounds of $a_{j,t}$, $\check{p}_{i,t}$, $\hat{p}_{i,t}$ and ${\beta}_{i,t}$, respectively. 
\end{theorem}

The proof of Theorem \ref{theorem_2} can be found in Appendix \ref{app_B}. Theorem \ref{theorem_2} tells us a large $V$ improves the performance of the proposed algorithm by reducing the optimality gap. However, the queue state can be very congested, potentially violating the carbon constraint. Conversely, a smaller $V$ results in a stable queue state but gives a larger optimality gap. Therefore, the parameter $V$ should be carefully selected.

\subsection{Disaggregation and queue update with dispatch signals}
In the last subsection, we quantify the power flexibility interval for each DC. The power grid operator is free to choose any power level within the reported flexibility interval $[\check{p}_{i,t},\hat{p}_{i,t}]$, and sends a dispatch signal $p_{i,t}^d$ to the DC operator to utilize this flexibility asset. The DC operator needs to determine the workload flow throughout the entire DC system to fulfill this dispatch signal (i.e., setting the power level of the back-end DC $i$ to be $p_{i,t}^d$), which we call ``disaggregation''. %We characterize the disaggregation feasibility of our proposed algorithm by the following proposition.
\begin{proposition}(Simultaneous one-slot disaggregation): The
Cartesian product of power intervals obtained from any feasible
solution of $\mathbf{P3}$  satisfies Definition  \ref{def-1}. \label{Prop-1}
\end{proposition}
\begin{proof}
Let
\(W_{i,t} := Q_{b,i,t}
+ \sum_{j\in\Omega_i^b}\check{m}_{ji,t}\)
and
\(H_{i,t} :=
\sum_{j\in\Omega_i^b}
(\hat{m}_{ji,t}-\check{m}_{ji,t})\).
Both are nonnegative. Define the routing interpolation factor
\begin{equation}
\lambda_{i,t}
=
\begin{cases}
\dfrac{\left[p_{i,t}^{d}\Delta t-W_{i,t}\right]^+}{H_{i,t}},
& H_{i,t}>0, \\[6pt]
0, & H_{i,t}=0.
\end{cases}
\tag{33}
\label{eq:routing-interpolation}
\end{equation}
Since
\(p_{i,t}^{d}\Delta t
\le \hat{p}_{i,t}\Delta t
\le W_{i,t}+H_{i,t}\),
we have \(0\le\lambda_{i,t}\le 1\).
if \(H_{i,t}=0\), the required processing is already
supported by \(W_{i,t}\). Choose
\begin{equation}
m_{ji,t}^{d}
=
(1-\lambda_{i,t})\check{m}_{ji,t}
+\lambda_{i,t}\hat{m}_{ji,t}.
\tag{34}
\label{eq:realized-routing}
\end{equation}
Thus
\(\check{m}_{ji,t}
\le m_{ji,t}^{d}
\le \hat{m}_{ji,t}
\le M_{ji}\),
and, crucially, for every front-end MN \(j\),
\(\sum_i m_{ji,t}^{d}
\le \sum_i\hat{m}_{ji,t}
\le Q_{f,j,t}+a_{j,t}\).
The independent factors \(\lambda_{i,t}\) therefore do not
overallocate shared front-end workload.
The resulting back-end workload satisfies
\begin{equation}
Q_{b,i,t}
+\sum_{j\in\Omega_i^b}m_{ji,t}^{d}
=
\max\left\{W_{i,t},\,p_{i,t}^{d}\Delta t\right\}.
\tag{35}
\label{eq:available-backend-workload}
\end{equation}
Consequently,
\begin{equation}
Q_{b,i,t}
+\sum_{j\in\Omega_i^b}m_{ji,t}^{d}
-p_{i,t}^{d}\Delta t
=
\left[W_{i,t}-p_{i,t}^{d}\Delta t\right]^+
\ge 0.
\tag{36}
\label{eq:nonnegative-backend-backlog}
\end{equation}
Also \(0\le p_{i,t}^{d}\le\overline{P}_i\).
All modeled one-slot constraints are therefore satisfied
simultaneously. %However, it is worth noting that this does not assert an unconditional delay or intertemporal carbon guarantee for arbitrary future dispatch sequences.
\end{proof}

To avoid the underestimation of future DC Power flexibility, when the dispatch signal is revealed, the actual queue states are updated with realized routing and power. For the workload queues, we have:
\begin{equation}\label{eq:update_1}
Q_{f,j,t+1} = Q_{f,j,t} + a_{j,t} - \sum\nolimits_{i \in \Omega_{j}^f} m_{ji,t}^d
\end{equation}
and
\begin{equation}\label{eq:update_2}
Q_{b,i,t+1} = Q_{b,i,t} + \sum\nolimits_{j \in \Omega_{i}^b} {m}_{ji,t}^d - {p}_{i,t}^d.
\end{equation}

For the virtual delay queues, we have:
\begin{equation}\label{eq:update_3}
Z_{f,j,t+1} = \max\{Z_{f,j,t} - \sum\nolimits_{i \in \Omega_{j}^f} {m}_{ji,t}^d + \delta_{f,j} \mathbbm{1}_{Q_{f,j,t}},0 \},
\end{equation}
and
\begin{equation}\label{eq:update_4}
Z_{b,i,t+1} = \max\{Z_{b,i,t} - {p}_{i,t}^d + \delta_{b,i} \mathbbm{1}_{Q_{b,i,t}},0 \}.
\end{equation}

For the virtual carbon queue, we have:
\begin{equation}\label{eq:update_5}
    Q_{c,t+1} = \max\{Q_{c,t} + \sum\nolimits_{i \in \mathcal{I}} \beta_{i,t} {p}_{i,t}^d - E, 0\}.
\end{equation}

With the above updated states of all queues, the DC operator can determine the flexibility interval when entering into the next time slot. In Algorithm \ref{alg_1}, we summarize the overall framework for power flexibility quantification of DCs.

\begin{algorithm}
\caption{Power Flexibility Quantification of DCs}
\label{alg_1}
\begin{algorithmic}
    % \Statex \hspace*{-2em} \textbf {\textit{III. Aggregation}}
    \State \textbf{Initialization:}
    \State $t \gets 1$. The DC operator initializes $Q_{f,j,1}$, $Q_{b,i,1}$,$Z_{f,j,1}$, $Z_{b,i,1}$, and $Q_{c,1}$ and sets the parameter $V$ and $\alpha$.
    \State \textbf{1. Flexibility Quantification:} 
    \State - Each MN collects workload requests $a_{j,t}$ for processing.
    % and append them to the queue $Q_{f,j,t}$.
     \State - The DC operator observes real-time gird carbon intensity of each back-end DC $\beta_{i,t}$, receives the reward signal $\gamma_{i,t}$ from each local grid operator, and solves $\mathbf{P3}$ to get the flexibility interval $ [\check{p}_{i,t}^\ast,\hat{p}_{i,t}^\ast]$ of each back-end DC.
    \State \textbf{2. Flexibility Provision:}
    \State - The DC operator reports $ [\check{p}_{i,t}^\ast,\hat{p}_{i,t}^\ast]$ to the corresponding grid operator.
    \State - The grid operator receives the flexibility interval and issues the dispatch signal $p_{i,t}^d \in[\check{p}_{i,t}^\ast,\hat{p}_{i,t}^\ast]$.
    \State \textbf{3. Disaggregation and  Queue Update:}
    \State - The DC operator performs the disaggregation process to determine $m_{ji,t}^d$ according to \eqref{eq:routing-interpolation} and \eqref{eq:realized-routing}.
    \State - The DC operator updates $Q_{f,j,t+1}$, $Q_{b,i,t+1}$, $Z_{f,j,t+1}$, and $Z_{b,i,t+1}$ {according to \eqref{eq:update_1}-\eqref{eq:update_4}}, and updates $Q_{c,t+1}$ according to \eqref{eq:update_5}.
    \State - Let $t \gets t+1$ and go back to step \textbf{1}.
\end{algorithmic}
\end{algorithm}

\section{Numerical simulation}
\label{sec:case}
In this section, we present numerical results to verify our proposed model and algorithm. 

\subsection{Simulation setup}
\begin{figure}[!htpb]
    \centering
    \includegraphics[width = 0.95\linewidth]{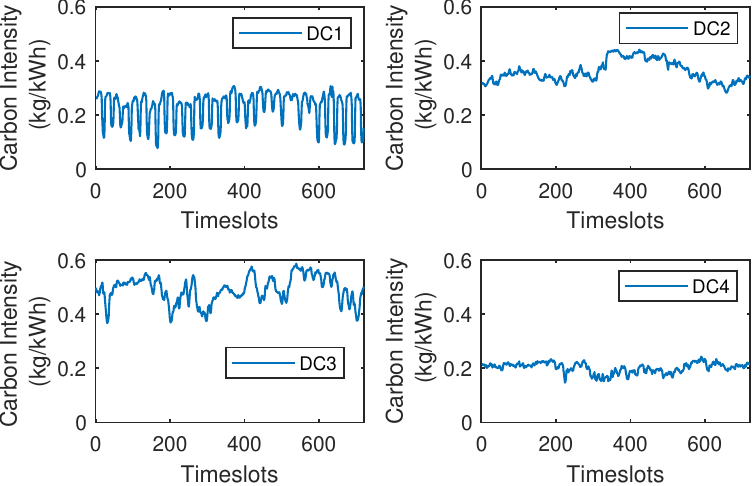} 
    % \vspace{-1em}
    \caption{ Carbon intensity profile of back-end DCs.}
    \label{fig:carbon_data}
\end{figure}

We conduct simulations over a period of one month with 720 one-hour slots.
We consider a DC system with 8 front-end MNs and 4 bank-end DCs and use real-world workload data of Google cluster from~\cite{google_workload}. %In Figure~\ref{fig:workload_data}, we display the workload profile of a front-end MN. 
The carbon intensity data of different electric regions is obtained from Electricity Maps~\cite{carbon_data}.
Figure~\ref{fig:carbon_data} shows the carbon intensity profiles of 4 back-end DCs, which varies significant across different locations.
Let $V$ = 15, and $\alpha$ = 1. The target carbon emission rate is $E$ = 10 t CO$_2$/h. 
For the DC system, let $M_{ji}$ = 2 MWh, $\overline{P}_i$ = 30 MWh, $\delta_{f,j}$ = 1 MWh, and $\delta_{b,i}$ = 2 MWh. 
The experiments are performed on the MATLAB 2024a platform with an Intel(R) Core(TM) Ultra 5 225H  processor and 32 GB RAM. 

We compare our approach with the following two online benchmarks. The detailed description of the two benchmark algorithms are given in Appendix \ref{app_C}.

\textbf{Static Myopic Algorithm (SMA)}: In this method, we equally allocate the total carbon budget to each time slot. Thus, the carbon emission limit of each time is $E$ since $\Delta t=1$ h. Under this carbon constraint, we greedily maximize the DC power flexibility value at each time slot with the objective function $\sum_{i \in \mathcal{I}} \gamma_{i,t}(\Hat{p}_{i,t} -\Check{p}_{i,t})$.  To maintain the queue stability of $Q_{f,j,t}$ and $Q_{b,i,t}$, queue upper limits are imposed on $Q_{f,j,t}$ and $Q_{b,i,t}$. For a fair evaluation, these bounds are selected to match the residual queue backlogs accumulated by our proposed algorithm at the end of the execution.

\textbf{Adaptive Myopic Algorithm (AMA)}: This approach is similar to SMA with a different strategy to handle the carbon emission constraint. Specifically, this approach utilizes the dispatch signals $p_{i,t}^d$ to adaptively update the time-average carbon emission limit. This can redistribute the previous unspent carbon budget for subsequent time slots.
% \begin{figure}[!htpb]
%     \centering
%     \includegraphics[width = 0.85\linewidth]{figure/workload.pdf} 
%     % \vspace{-1em}
%     \caption{ Received workload profile of a front-end node.}
%     \label{fig:workload_data}
% \end{figure}

\subsection{Flexibility quantification results}
\begin{figure}[!htpb]
    \centering
    \includegraphics[width = 0.95\linewidth]{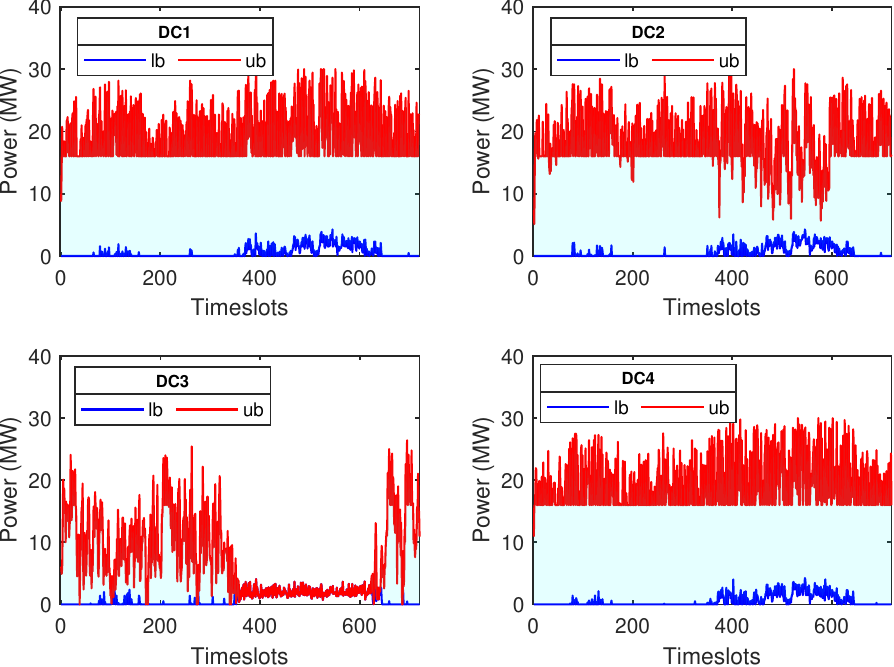} 
    % \vspace{-1em}
    \caption{ Flexibility intervals of back-end DCs.}
    \label{fig:flexibility_regions}
\end{figure}

\begin{figure}[!htpb]
    \centering
    \includegraphics[width = 0.8\linewidth]{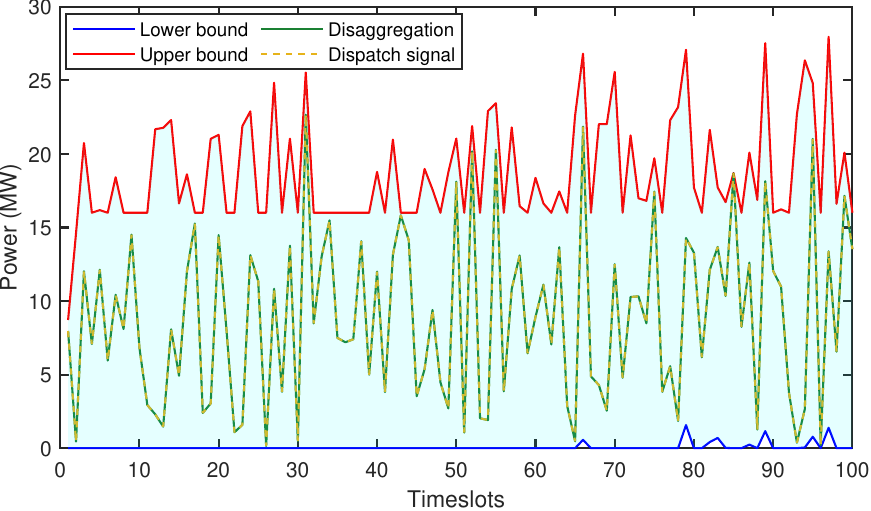} 
    % \vspace{-1em}
    \caption{ Disaggregation result of a back-end node.}
    \label{fig:dispatch_signal}
\end{figure}

% Note that the dispatch signals are normalized to the range of [0,1] based on their real-time lower and upper bounds $\check{p}_{i,t}$ and $\hat{p}_{i,t}$. 
We run the proposed algorithm to derive the flexibility intervals with the corresponding dispatch signals randomly generated. We first consider the setting when $\gamma_{i,t}=1,\forall i,t$.
%We randomly generate the dispatch signals to activate the flexibility quantification process. Then 
The flexibility intervals of 4 back-end DCs are shown in Figure~\ref{fig:flexibility_regions}. The flexibility intervals are defined by two trajectories: the upper bound (red line) and the lower bound (blue line). The space enclosed between the upper and lower bounds is the flexibility interval for each time slot. We can observe that the lower bounds consistently maintain low values below 5 MW, even dropping to zero for certain periods. In contrast, the upper bounds are notably higher and highly fluctuating. These phenomena can be interpreted with reference to the dispatch signal displayed in Figure~\ref{fig:dispatch_signal} (yellow dashed line). To be specific, the dispatch signals issued by grid operators determine the actual power consumption of DCs, and consequently governs the resulting carbon emissions. Owing to the high volatility of dispatch signals, the overall system's carbon emissions exhibit substantial fluctuations over time. To control the time-average emissions, the proposed algorithm adaptively constrain the upper bounds, causing them to fluctuate accordingly. {The pattern of lower bounds can be explained by the fact that higher power levels are likely to be frequently dispatched, as indicated by the dispatch signals.} This actually promotes the timely processing of workloads, thereby preventing excessive delays. Consequently, the lower bounds can be maintained at a low level to better utilize the deferrable nature of workloads. Moreover, Figure~\ref{fig:dispatch_signal} also presents the disaggregation results (depicted by the green line). It can be observed that the disaggregation curve completely coincides with the curve of dispatch signals, which demonstrates that our algorithm ensures the feasibility of disaggregation and thereby validates Proposition~\ref{def-1}.

TABLE~\ref{tab:comp_DCs} displays both the cumulative flexibility value of each back-end DC and its regional mean carbon intensity. It is clear that the DC, which belongs to the region with a lower mean carbon intensity, achieves a higher level of cumulative flexibility. This observation can be attributed to the fact that under the system-level emission constraint, the proposed algorithm achieves flexibility maximization through raising the upper bounds of the flexibility intervals for DCs located in regions with low carbon intensity. In this way, workloads are more likely to be transferred to clean-energy electric regions.

\begin{table}[t]
\small
    \centering
    \caption{Comparison of cumulative flexibility and regional mean carbon intensity across different back-end DCs.} 
    % \vspace{-0.5em}
    \begin{tabular}{@{}ccccccccc@{}}
        \toprule
       \multirow[c]{2}{*}{} & Cumulative  & Mean carbon  \\
       & flexibility value (MWh) & intensity (kg CO$_2$/kWh)\\
        \midrule
        DC1  & 13645 & 0.2282 \\
        DC2  & 12694 & 0.3607 \\
        DC3 & 4690& 0.4911\\
        DC4 &13672 &0.2016\\
        \bottomrule
        \label{tab:comp_DCs}
   \end{tabular}
\end{table}

\begin{figure}[!htpb]
    \centering
    \includegraphics[width = 0.8\linewidth]{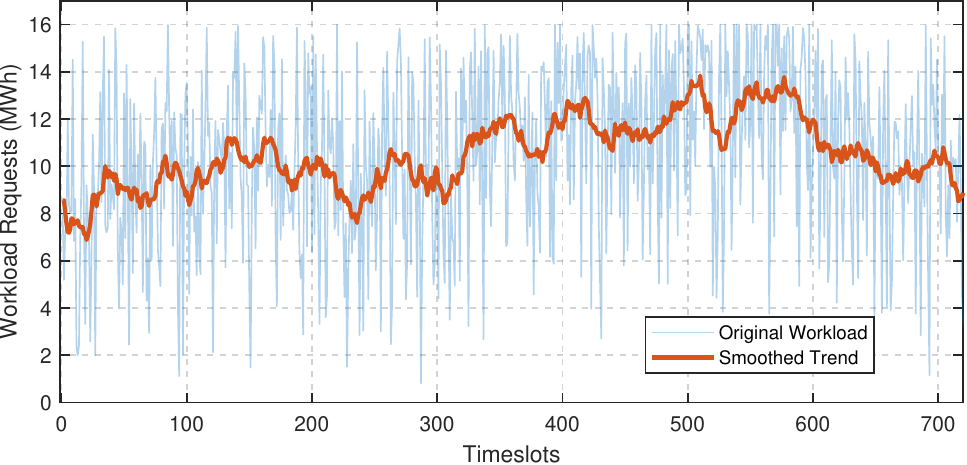} 
    % \vspace{-1em}
    \caption{ Received workload profile of a back-end node.}
    \label{fig:workload_back}
\end{figure}

\begin{figure}[!htpb]
    \centering
    \includegraphics[width = 0.95\linewidth]{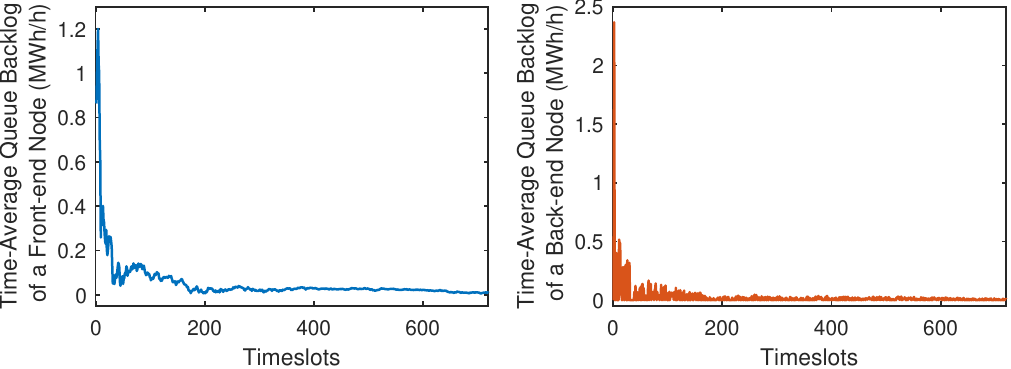} 
    % \vspace{-1em}
    \caption{ Queue backlogs of a front-end node and a back-end node.}
    \label{fig:queue_backlog}
\end{figure}

One may notice that in Figure~\ref{fig:flexibility_regions}, there seems to be a clear lower bound (i.e., 16 MW) beneath the oscillating flexibility upper bounds. This phenomenon can be explained by Figure~\ref{fig:workload_back}. In Figure~\ref{fig:workload_back}, we display the received workload of a back-end DC with its smoothed trend. As our DC system has 8 front-end nodes, the maximum workload that can be transferred to a back-end node each time is 16 MWh under the link bandwidth of 2 MWh. Thus, once the previous stored workloads of back-end DCs are fully depleted, the attainable flexibility upper bound becomes strictly constrained by the maximum workload received from the front-end nodes. In Figure~\ref{fig:queue_backlog}, we also show the time-average queue backlogs of a front-end node and a back-end node over time. We see that the time-average queue backlogs of both the front-end and back-end nodes eventually stabilize near zero, which demonstrates the superior performance of our proposed algorithm in maintaining queue stability.

% When the $\gamma_{i,t}$ is set to be equal for all back-end DCs, the cumulative flexibility across different DCs is determined by regional carbon intensity level. To evaluate the impact of $\gamma_{i,t}$, we set different values for  different DCs to distinguish the flexibility value of each DC. The ratio of the $\gamma_{i,t}$ is 1:2:3:4 for DC1, DC2, DC3, and DC4, respectively. Figure~\ref{fig:flexibility_regions_1} shows the corresponding results. One can observe that the distribution of flexibility is no longer determined by regional carbon intensity, but by the value of $\gamma_{i,t}$. Although the region where DC3 is located has the highest carbon intensity level, its cumulative flexibility  exceeds those of DC1 and DC2 due to the higher reward signal $\gamma_{i,t}$ assigned to DC3. Thus, the spatio-temporal distribution of the flexibility within the DC system can be regulated by setting the value of $\gamma_{i,t}$.
% However, since we place greater emphasis on the flexibility of high carbon intensity regions (i.e., DC3), the cumulative flexibility of the entire DC system  decreases under fixed carbon emission budget.
% % total flexibility 32805
% \begin{figure}
%     \centering
%     \includegraphics[width = 0.95\linewidth]{figure/flexibility_regions_1.pdf} 
%     % \vspace{-1em}
%     \caption{ Flexibility intervals of back-end DCs .}
%     \label{fig:flexibility_regions_1}
% \end{figure}

\begin{figure}[!htpb]
    \centering
    \includegraphics[width = 0.75\linewidth]{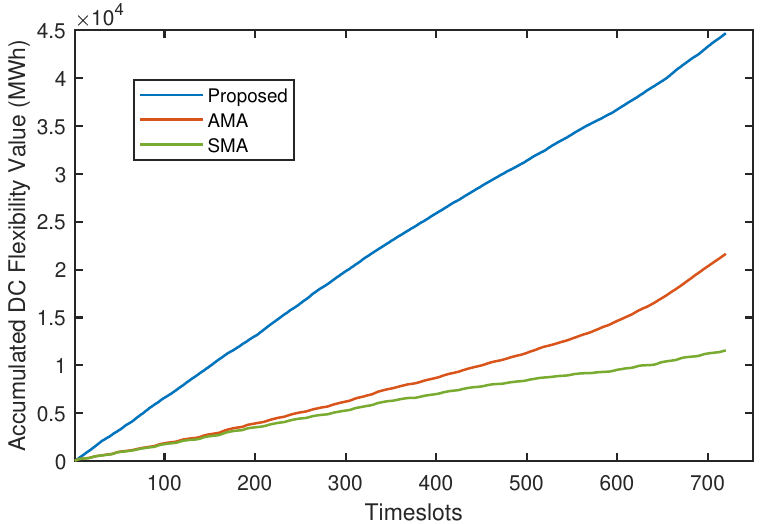} 
    % \vspace{-1em}
    \caption{Time accumulated DC flexibility vaule by different methods when the reward signals are set to a constant of 1.}
    \label{fig:flex_area_cmp}
\end{figure}

Figure~\ref{fig:flex_area_cmp} compares the accumulated DC flexibility value of our method against two baselines. Our approach consistently outperforms both baselines, as it not only leverages the deferrable nature of workloads to control the lower bound, but also efficiently utilizes the limited carbon budget for the upper bound control under time-varying carbon intensity. Both our method and AMA outperform SMA, as both leverage dispatch signals to update carbon constraints to fully utilize the limited carbon budget.
TABLE~\ref{tab:flex_area_method_comp}  further compares the total flexibility value, time-average emission rate, and performance index across the three methods, where the performance index is defined as the ratio of total flexibility value to total emissions. It is clear that all three methods satisfy the carbon constraints. However, AMA and SMA resort to myopic strategies to avoid violations of carbon constraints, rendering them incapable of fully utilizing the available carbon budget to improve flexibility. Moreover, our method achieves a markedly higher performance index than both baselines, indicating its superior ability to leverage the carbon budget under complex and time-varying environments to fully exploit the flexibility of DCs.

\begin{figure}[!htpb]
    \centering
    \includegraphics[width = 0.75\linewidth]{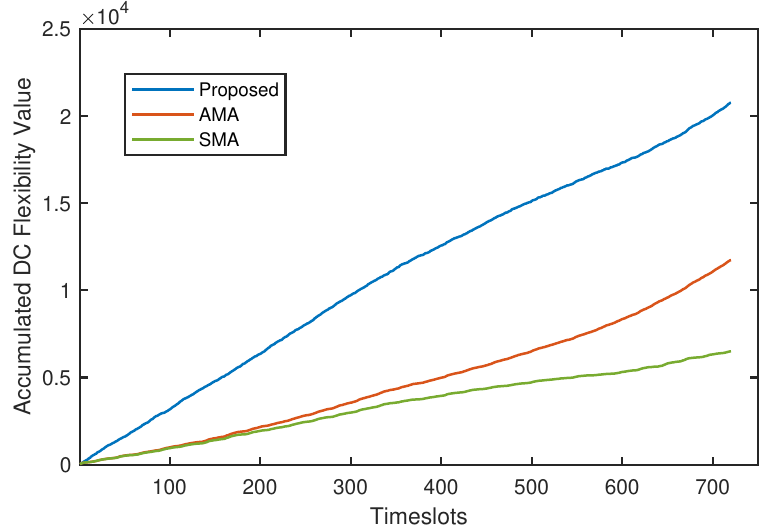} 
    % \vspace{-1em}
    \caption{Time accumulated DC flexibility value by different methods when the rewards signals vary spatiotemporally.}
    \label{fig:flex_value_cmp}
\end{figure}

We next consider a scenario where reward signals (i.e.,$\gamma_{i,t},\forall i,t$) vary spatiotemporally, meaning that the value of flexibility differs across time and electric regions. For simplicity, reward signals for different time slots and regions are randomly generated within $[0, 1]$. Here, the unit for this signal is per MWh. The corresponding comparison results are presented in Figure~\ref{fig:flex_value_cmp} and TABLE~\ref{tab:flex_value_method_comp}. One can observe that the proposed algorithm effectively responds to the time-varying reward signals, achieving a higher flexibility value than both baselines. Furthermore, the performance index reported in TABLE~\ref{tab:flex_value_method_comp} further corroborate that our approach can unlock the value of DC flexibility under various uncertainties.

% the cumulative carbon emissions of the proposed method with and without emission control. One can observe 
% an emission reduction of around 500 t is achieved through the established virtual carbon queue, compared to the case without emission control.

% \begin{figure}
%     \centering
%     \includegraphics[width = 0.85\linewidth]{figure/emission_control_cmp.pdf} 
%     % \vspace{-1em}
%     \caption{Cumulative emissions with and without carbon emission control.}
%     \label{fig:emission_cmp}
% \end{figure}

\begin{table}[t]
\small
    \centering
    \caption{Results of DC flexibility value
    under different methods when the reward signals are set to a constant of 1.} 
    \vspace{-0.5em}
    \begin{tabular}{@{}ccccccccc@{}}
        \toprule
        Method & Proposed & SMA & AMA \\
        \midrule
        Total flexibility value (MWh) & 44701 & 21652  & 11557 \\
        Carbon emission rate (t CO$_2$/h)  & 10.00 & 9.14 & 8.39 \\
        % Unfulfilled energy (kWh)  & 0 & 218.6 & 2.2 & 10.0\\ 
        % \midrule
        Performance index & 6.21 & 3.29 & 1.91 \\
        \bottomrule
    \label{tab:flex_area_method_comp}
   \end{tabular}
\end{table}

\begin{table}[t]
\small
    \centering
    \caption{Results of DC flexibility value under different methods when the rewards signals vary spatiotemporally} 
    \vspace{-0.5em}
    \begin{tabular}{@{}ccccccccc@{}}
        \toprule
        Method & Proposed & SMA & AMA \\
        \midrule
        Total flexibility value & 20759 & 11735 & 6504.7\\
        Carbon emission rate (t CO$_2$/h)  & 9.99 & 9.21 & 8.54 \\
        % Unfulfilled energy (kWh)  & 0 & 218.6 & 2.2 & 10.0\\ 
        % \midrule
        Performance index & 2.89 & 1.77 &1.06 \\
        \bottomrule
    \label{tab:flex_value_method_comp}
   \end{tabular}
\end{table}

\subsection{Impact of factors}
\begin{figure}[!htpb]
    \centering
    \includegraphics[width = 0.75\linewidth]{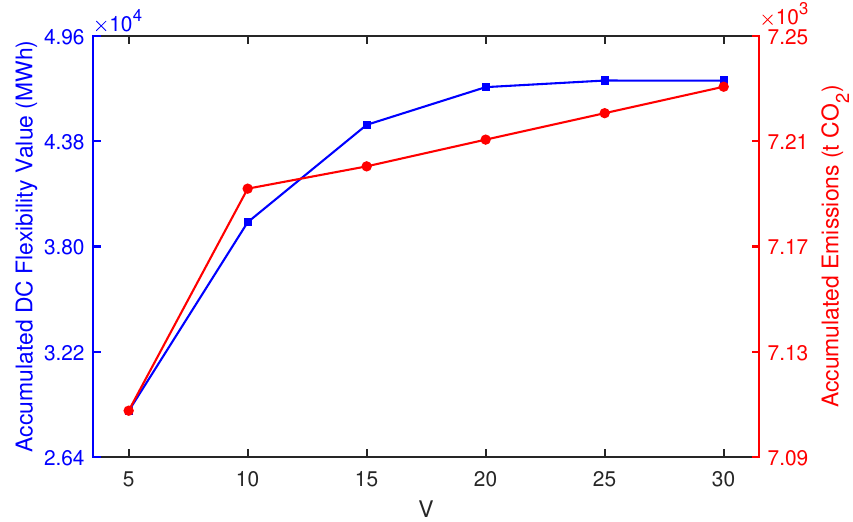} 
    % \vspace{-1em}
    \caption{ Impact of parameter $V$.}
    \label{fig:impact_V}
\end{figure}

The impact of several factors including parameter $V$, $\alpha$, and $E$ is tested in this subsection. We vary the value of $V$ in a range from 5 to 30 to evaluate its effect, and the results are shown in Figure~\ref{fig:impact_V}. We see that both cumulative flexibility value and emissions exhibit a monotonic increase with rising $V$, which suggests a trade-off between flexibility provision and carbon emissions. This trend stems from the role of $V$, which is to adjust the significance of flexibility maximization in the objective function of $\mathbf{P3}$. This insight is consistent with the analysis in Theorem \ref{theorem_2}.

\begin{figure}[!htpb]
    \centering
    \includegraphics[width = 0.75\linewidth]{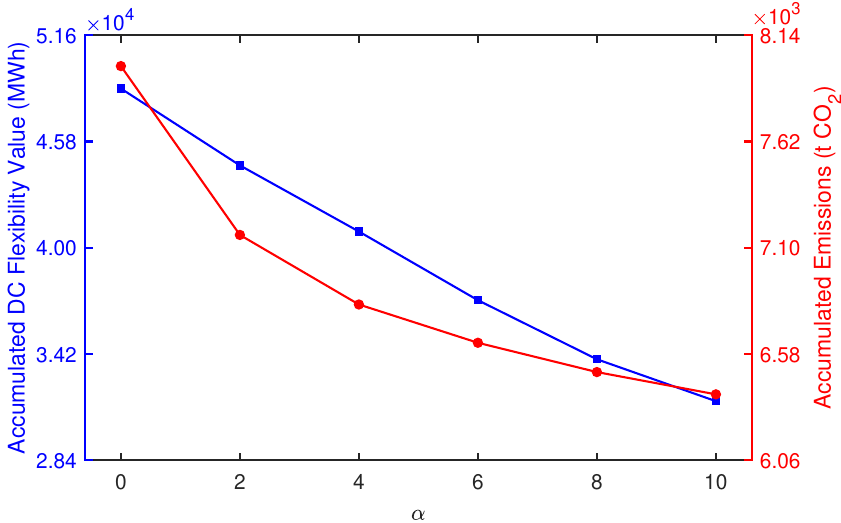} 
    % \vspace{-1em}
    \caption{ Impact of parameter $\alpha$.}
    \label{fig:impact_alpha}
\end{figure}

A sensitive analysis is conducted by changing the value of $\alpha$, and the results are presented in Figure~\ref{fig:impact_alpha}. It can be seen that cumulative flexibility value and emissions decrease with increasing $\alpha$, which again implies the trade-off between flexibility provision and carbon emissions. The underlying reason is that when increasing $\alpha$, the relative importance of the carbon queue in the Lyapunov function \eqref{eq:lay_func} is raised. Therefore,  the algorithm performs stronger stability control over the carbon queue than over the workload queues, which is responsible for the pattern observed in Figure~\ref{fig:impact_alpha}.

\begin{figure}[!htpb]
    \centering
    \includegraphics[width = 0.75\linewidth]{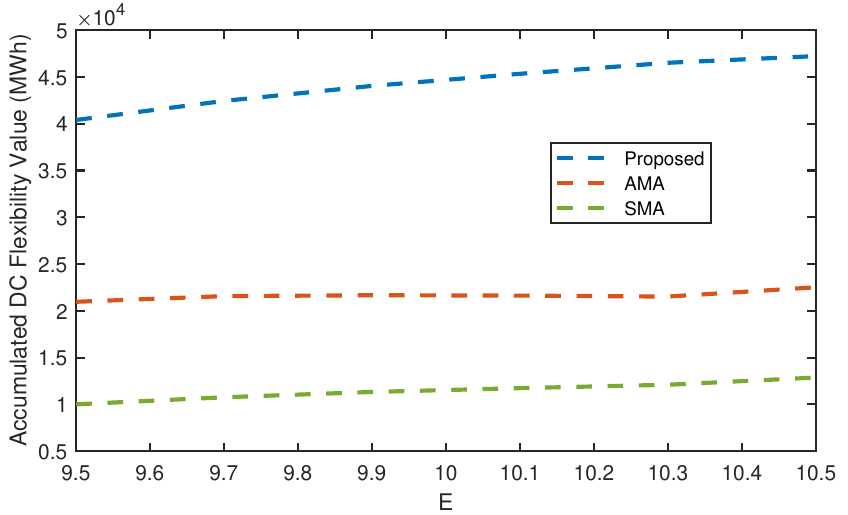} 
    % \vspace{-1em}
    \caption{ Time accumulated DC flexibility value by different methods under varying parameter $E$.}
\label{fig:accu_flex_area_diff_carbon}
\end{figure}

We assess the performance of the algorithm under different carbon emission rate constraints by changing $E$. Figure~\ref{fig:accu_flex_area_diff_carbon} shows that as $E$ increases, the cumulative flexibility value of our proposed method rises. This is straightforward since a larger $E$ imposes a less strict carbon constraint on the entire DC system, so that the DCs are able to explore more flexibility under more carbon budget. Moreover, our approach remains consistently superior to AMA and SMA in terms of DC flexibility value over the variation range of $E$. AMA and SMA are relatively insensitive to variations in $E$ due to their myopic strategies that fail to fully explore the DC power flexibility under varying parameter $E$.

\begin{figure}[!htpb]
    \centering
    \includegraphics[width = 0.75\linewidth]{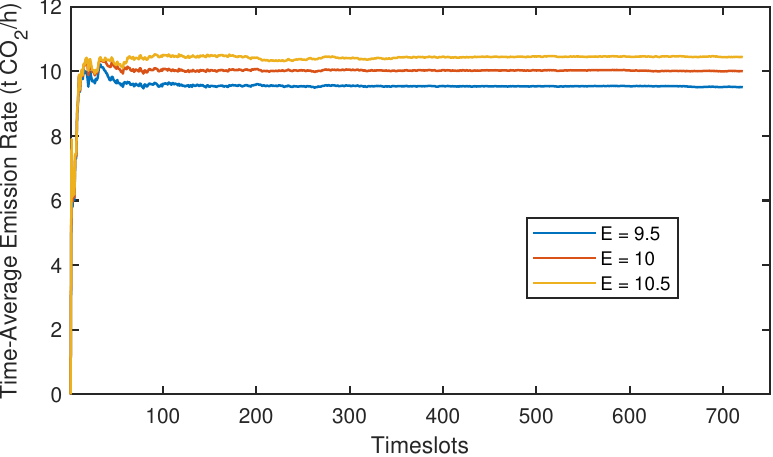} 
    % \vspace{-1em}
    \caption{ Time-average carbon emission rate of the proposed method under different parameter $E$.}
    \label{fig:emission_profiles}
\end{figure}

In TABLE~\ref{tab:impact_carbon}, we further show the results of our proposed method under different parameter $E$, where emission rate refers to the time-average emission rate of the DC system over all timeslots. We see that our approach effectively adheres to the prescribed carbon emission rate limit $E$. We also illustrate the running average of the proposed algorithm's carbon emission rate over time in Figure~\ref{fig:emission_profiles}. One may observe that our algorithm effectively brings the time-average carbon emission rate into compliance with the carbon constraint, and maintains its long-term stability throughout the whole time horizon.

\begin{table}[t]
\small
    \centering
    \caption{Cumulative flexibility and emissions of proposed method under different carbon constraints.} 
    % \vspace{-0.5em}
    \begin{tabular}{@{}ccccccccc@{}}
        \toprule
       \multirow{2}{*}{E } & Cumulative  & Cumulative & Emission\\
        & flexibility & emissions (t CO$_2$) & rate (t CO$_2$/h)\\
        \midrule
        9.5  & 40398 & 6845.3 & 9.5 \\
        10   & 44701 & 7200.9 & 10.0\\
        10.5 & 47223 & 7515.9 & 10.4\\
        \bottomrule
        \label{tab:impact_carbon}
   \end{tabular}
\end{table}

\subsection{Scalability analysis}
We perform a scalability analysis to demonstrate the computational efficiency of the proposed method. We evaluate the average runtime per time slot for both flexibility quantification and disaggregation steps under different numbers of MNs and DCs, where we increase the DC system scale with a fixed MN-to-DC ratio of 2:1. The corresponding results are summarized in TABLE~\ref{tab:sca_ana}. As observed, the runtime of both steps increases with the system scale. For flexibility quantification, even when the problem scale reaches a total number of 96 for MNs and DCs, its execution time remains slightly over 2 seconds, which is well below the time interval (i.e., 1 hour). For disaggregation step, the execution time is on the order of $10^{-5}$ seconds, making it practically negligible. Consequently, these results demonstrate the superior computational efficiency of our proposed algorithm.
\begin{table}[t]
\small
    \centering
    \caption{Runtimes of the proposed algorithm under different numbers of MNs and DCs} 
    \vspace{-0.5em}
    \begin{tabular}{@{}ccccccccc@{}}
        \toprule
        Total number of MNs and DCs &12 & 24 & 48 & 96 \\
        \midrule
        Flexibility quantification (s) & 0.031 & 0.046 & 0.21 & 2.66\\
         Disaggregation ($\times 10^{-5}$ s)  & 1.40 & 1.70 & 2.84&5.48 \\
        \bottomrule
    \label{tab:sca_ana}
   \end{tabular}
\end{table}

\section{Conclusion}
\label{sec:conclusion}
In this paper, we quantify the power flexibility of DCs in real time under carbon emission constraints. We adopt flexibility intervals to describe the flexibility potential of DCs. A set of virtual queues are constructed to control the DC workload processing delay and carbon emission rate.
Then, we develop an online algorithm to calculate the flexibility intervals, considering uncertainties in both workload arrival and grid carbon intensity. Our algorithm is prediction-free and independent of modeling of uncertainty. The efficiency of our online method are validated through simulation studies.
Specifically, our algorithm achieves an effective flexibility quantification without violating the carbon constraints while outperforming benchmark algorithms.

In future work, we will extend our online power flexibility quantification framework to incorporate more detailed workload characteristics.

%other loads with deferrable nature.} Another possible direction is to consider the participation of DCs in a flexibility market.

% Appendix
%-------------------------------------------
% Activate the appendix in the doc
% from here on sections are numerated with capital letters 
%\appendix

% Change equation numbering format to be sequential within sections in the appendix
\renewcommand\theequation{\Alph{section}.\arabic{equation}} % Redefine equation numbering format
\counterwithin*{equation}{section} % Number equations within sections
\renewcommand\thefigure{\Alph{section}\arabic{figure}} % Redefine equation numbering format
\counterwithin*{figure}{section} % Number equations within sections
\renewcommand\thetable{\Alph{section}\arabic{table}} % Redefine equation numbering format
\counterwithin*{table}{section} % Number equations within sections

% % \newpage

\begin{appendices}

\section{Proof of Theorem \ref{theorem_1}}\label{app_A}
We first prove that a workload request that arrives at back end $i$ will wait for at most $D_{i,b}$ before being finished. Without loss of generality, we suppose a request arrives at back end $i$ at time $t_0$, and it will be finished exactly at time $t_0 + \tau$. Therefore, within the time period [$t_0$, $t + \tau-1$], this request has not been completed, which indicates that the queue backlog $Q_{i,b,t} > 0 $ holds during the time period [$t_0$, $t_0 + \tau-1$]. Thus, we have:
\begin{equation}
    Z_{b,i,t+1} \geq {Z_{b,i,t}} -\check{p}_{i,t} + \delta_{b,i} , \forall t \in [t_0,t_0+\tau-1]
\end{equation}

By summing the above over $t \in [t_0,t_0+\tau-1]$, we have:

\begin{equation}
    Z_{b,i,t_0+\tau} - Z_{b,i,t_0} \geq -\sum_{t = t_0}^{t_0+\tau-1} \check{p}_{i,t} + \tau \delta_{b,i}
\end{equation}

Since $0 \leq Z_{b,i,t} \leq Z_{b,i}^{max}$, we have:
\begin{equation}\label{eq:A_3}
    \sum_{t = t_0}^{t_0+\tau-1} \check{p}_{i,t} \geq \tau \delta_{b,i} - Z_{b,i}^{max}
\end{equation}

The request is appended to the queue backlog $Q_{i,b,t}$ at time $t$ and is not completely cleared until time $t + \tau$. Since the request is processed in FIFO manner, the amount of the requests that DC $i$ finishes during the time period [$t_0$, $t_0 + \tau-1$ ] must be less than the upper bound of $Q_{b,i,t}$, i.e.,

\begin{equation}\label{eq:A_4}
    \sum_{t = t_0}^{t_0+\tau-1} \check{p}_{i,t} < Q_{b,i}^{max}
\end{equation}

Combining \eqref{eq:A_3} and \eqref{eq:A_4} yields:

\begin{equation}\label{eq:A_5}
    \tau < (Q_{b,i}^{max} + Z_{b,i}^{max})/\delta_{b,i}
\end{equation}

The above equation shows that the worst-case request delay at back end $i$ is $D_{b,i}$. With a similar procedure, we can prove that the largest request delay at the front end $j$ is $D_{f,j}$, where $D_{f,j} = (Q_{f,j}^{max} + Z_{f,j}^{max})/\delta_{f,j}$. As a result, the worst-case request delay of front end $j$ is $D_{f,j} + \max_{i \in \Omega_j^f}D_{b,i}$.

\section{Proof of Theorem \ref{theorem_2}}\label{app_B}
According to the Lyapunov optimization theory~\cite{neely2010stochastic},
for any given $\Tilde{\varepsilon} > 0 $, there exists a stationary $w$-only strategy $\check{p}_{i,t}^w, \hat{p}_{i,t}^w, \forall i,t$, which is feasible to $\mathbf{P2}$ with performance guarantee $\Tilde{\varepsilon}$. We denote the corresponding time-average power flexibility expectation by $\Tilde{f}$. The $w$-only strategy is independent of the information of queues and satisfies $f_1^\ast \geq \Tilde{f} \geq f_1^\ast - \Tilde{\varepsilon}$.

Given that $\check{p}_{i,t}^\ast, \hat{p}_{i,t}^\ast, \forall i,t$ is the optimal solution of $\mathbf{P3}$, we can have the following:
\begin{align*}
\label{eq_B1}
         &\Delta (\mathbf{{\Theta}}_t^\ast) +  V\sum\nolimits_{i \in \mathcal{I}} \gamma_{i,t}(\check{p}_{i,t}^\ast -\hat{p}_{i,t}^\ast) \leq  V\sum\nolimits_{i \in \mathcal{I}} \gamma_{i,t} (\check{p}_{i,t}^w -\hat{p}_{i,t}^w)\\
         & +\frac{1}{2}\sum_{j \in \mathcal{J}}\bigg[(a_{j,t} - \sum\nolimits_{i \in \Omega_j^f}\check{m}_{ji,t}^w)^2 +  ( \delta_{f,j}- \sum\nolimits_{i \in \Omega_j^f}\check{m}_{ji,t}^w)^2 \bigg] \\
         & +\frac{1}{2}\sum_{i \in \mathcal{I}}\bigg[(\sum\nolimits_{j \in \Omega_i^b}\check{m}_{ji,t}^w - \check{p}_{i,t}^w)^2 +  ( \delta_{b,i}- \check{p}_{i,t}^w)^2 \bigg] \\
         & +\sum_{j \in \mathcal{J}}\bigg[Q_{f,j,t}a_{j,t} + Z_{f,j,t}\delta_{f,j} - (Q_{f,j,t} + Z_{f,j,t})\sum_{i \in \Omega_j^f}\check{m}_{ji,t}^w \bigg]\\
         & + \sum_{i \in \mathcal{I}}\bigg[Q_{b,i,t}\sum_{j \in \Omega_i^b}\check{m}_{ji,t}^w + Z_{b,i,t}\delta_{b,i} - (Q_{b,i,t} + Z_{b,i,t})\check{p}_{i,t}^w \bigg]\\
         & + \frac{1}{2} \alpha ( \sum\nolimits_{i \in \mathcal{I}} \beta_{i,t}\hat{p}_{i,t}^w - E)^2 + \alpha Q_{c,t}(\sum\nolimits_{i \in \mathcal{I}} \beta_{i,t}\hat{p}_{i,t}^w - E)\\
         & \leq V\sum\nolimits_{i \in \mathcal{I}} \gamma_{i,t} (\check{p}_{i,t}^w -\hat{p}_{i,t}^w) \\
         & + \frac{1}{2}\sum_{j \in \mathcal{J}}\bigg[\max(A_{j}^2 , \sum\nolimits_{i \in \Omega_j^f}M_{ji}^2) + \max ( \delta_{f,j}^2, \sum\nolimits_{i \in \Omega_j^f}M_{ji}^2) \bigg] \\
         & + \frac{1}{2}\sum_{i \in \mathcal{I}}\bigg[\max(\sum\nolimits_{j \in \Omega_i^b}M_{ji}^2 , \check{P}_{i}^2) +  \max( \delta_{b,i}^2,
         \check{P}_{i}^2) \bigg]\\
         & + \sum_{j \in \mathcal{J}}\bigg[Q_{f,j,t}(a_{j,t}-\sum_{i \in \Omega_j^f}\check{m}_{ji,t}^w)\\
         &+ Z_{f,j}^{max}\delta_{f,j} - Z_{f,j,t}\sum_{i \in \Omega_j^f}\check{m}_{ji,t}^w \bigg] \stepcounter{equation}\tag{\theequation}\\ 
         & +\sum_{i \in \mathcal{I}}\bigg[Q_{b,i,t}(\sum_{j \in \Omega_i^b}\check{m}_{ji,t}^w - \check{p}_{i,t}^w) + Z_{b,i}^{max}\delta_{b,i} - Z_{b,i,t}\check{p}_{i,t}^w \bigg]\\
         & + \frac{1}{2} \alpha \max( \sum\nolimits_{i \in \mathcal{I}} \bar{\beta}_{i}^2\hat{P}_{i}^2, E^2) + \alpha Q_{c,t}(\sum\nolimits_{i \in \mathcal{I}} \beta_{i,t}\hat{p}_{i,t}^w - E)\\
          & \leq B + V\sum\nolimits_{i \in \mathcal{I}} \gamma_{i,t} (\check{p}_{i,t}^w -\hat{p}_{i,t}^w) + \sum_{j \in \mathcal{J}}Q_{f,j,t}(a_{j,t}-\sum_{i \in \Omega_j^f}\check{m}_{ji,t}^w)\\
          &  +\sum_{i \in \mathcal{I}}Q_{b,i,t}(\sum_{j \in \Omega_i^b}\check{m}_{ji,t}^w - \check{p}_{i,t}^w) + \alpha Q_{c,t}(\sum\nolimits_{i \in \mathcal{I}} \beta_{i,t}\hat{p}_{i,t}^w - E)
\end{align*}

For the $w$-only strategy, it satisfies the following:
\begin{equation}\label{eq:B2}
    \begin{aligned}
        & \mathbb{E}\big[(a_{j,t}-\sum\nolimits_{i \in \Omega_j^f}\check{m}_{ji,t}^w)\big] = 0 \\
        & \mathbb{E}\big[(\sum\nolimits_{j \in \Omega_i^b}\check{m}_{ji,t}^w - \check{p}_{i,t}^w)\big] = 0\\
        & \mathbb{E}\big[(\sum\nolimits_{i \in \mathcal{I}} \beta_{i,t}\hat{p}_{i,t}^w - E)\big] \leq 0 
    \end{aligned}
\end{equation}

Armed with \eqref{eq:B2}, we then take expectations on both sides of \eqref{eq_B1} and obtain the following by using the independence of the $w$-only strategy on queues:
\begin{align*}
\label{eq_B3}
 & \mathbb{E}\big[\Delta (\mathbf{{\Theta}}_t^\ast) +  V\sum\nolimits_{i \in \mathcal{I}} \gamma_{i,t}(\check{p}_{i,t}^\ast -\hat{p}_{i,t}^\ast)\big]\\
        & \leq B + V\mathbb{E}\big[\sum\nolimits_{i \in \mathcal{I}} \gamma_{i,t}(\check{p}_{i,t}^w -\hat{p}_{i,t}^w)\big]\\
        & +  \sum\nolimits_{j \in \mathcal{J}}\mathbb{E}\big[Q_{f,j,t}(a_{j,t}-\sum\nolimits_{j \in \Omega_i^b}\check{m}_{ji,t}^w)\big] \\
        & + \sum\nolimits_{i \in \mathcal{I}}\mathbb{E}\big[Q_{b,i,t}(\sum\nolimits_{j \in \Omega_i^b}\check{m}_{ji,t}^w - \check{p}_{i,t}^w)\big]\\
        & + \alpha \mathbb{E}\big[Q_{c,t}(\sum\nolimits_{i \in \mathcal{I}} \beta_{i,t}\hat{p}_{i,t}^w - E)\big] \\
        & = B + V\mathbb{E}\big[\sum\nolimits_{i \in \mathcal{I}} \gamma_{i,t}(\check{p}_{i,t}^w -\hat{p}_{i,t}^w)\big]\stepcounter{equation}\tag{\theequation}\\
        & +  \sum\nolimits_{j \in \mathcal{J}}\mathbb{E}\big[Q_{f,j,t}\big]\mathbb{E}\big[(a_{j,t}-\sum\nolimits_{j \in \Omega_i^b}\check{m}_{ji,t}^w)\big] \\
        & + \sum\nolimits_{i \in \mathcal{I}}\mathbb{E}\big[Q_{b,i,t}\big]\mathbb{E}\big[(\sum\nolimits_{j \in \Omega_i^b}\check{m}_{ji,t}^w - \check{p}_{i,t}^w)\big]\\
        & + \alpha \mathbb{E}\big[Q_{c,t}\big]\mathbb{E}\big[(\sum\nolimits_{i \in \mathcal{I}} \beta_{i,t}\hat{p}_{i,t}^w - E)\big] \\
        & \leq B + V\mathbb{E}\big[\sum\nolimits_{i \in \mathcal{I}} \gamma_{i,t}(\check{p}_{i,t}^w -\hat{p}_{i,t}^w)\big]
\end{align*}
% \begin{equation}
%     \begin{aligned}
%         & \mathbb{E}\big[\Delta (\mathbf{{\Theta}}_t^\ast) +  V\sum\nolimits_{i \in \mathcal{I}} \gamma_{i,t}(\check{p}_{i,t}^\ast -\hat{p}_{i,t}^\ast)\big]\\
%         & \leq B + V\mathbb{E}\big[\sum\nolimits_{i \in \mathcal{I}} \gamma_{i,t}(\check{p}_{i,t}^w -\hat{p}_{i,t}^w)\big]\\
%         & +  \sum\nolimits_{j \in \mathcal{J}}\mathbb{E}\big[Q_{f,j,t}(a_{j,t}-\sum\nolimits_{j \in \Omega_i^b}\check{m}_{ji,t}^w)\big] \\
%         & + \sum\nolimits_{i \in \mathcal{I}}\mathbb{E}\big[Q_{b,i,t}(\sum\nolimits_{j \in \Omega_i^b}\check{m}_{ji,t}^w - \check{p}_{i,t}^w)\big]\\
%         & + \alpha \mathbb{E}\big[Q_{c,t}(\sum\nolimits_{i \in \mathcal{I}} \beta_{i,t}\hat{p}_{i,t}^w - E)\big] \\
%         & = B + V\mathbb{E}\big[\sum\nolimits_{i \in \mathcal{I}} \gamma_{i,t}(\check{p}_{i,t}^w -\hat{p}_{i,t}^w)\big]\\
%         & +  \sum\nolimits_{j \in \mathcal{J}}\mathbb{E}\big[Q_{f,j,t}\big]\mathbb{E}\big[(a_{j,t}-\sum\nolimits_{j \in \Omega_i^b}\check{m}_{ji,t}^w)\big] \\
%         & + \sum\nolimits_{i \in \mathcal{I}}\mathbb{E}\big[Q_{b,i,t}\big]\mathbb{E}\big[(\sum\nolimits_{j \in \Omega_i^b}\check{m}_{ji,t}^w - \check{p}_{i,t}^w)\big]\\
%         & + \alpha \mathbb{E}\big[Q_{c,t}\big]\mathbb{E}\big[(\sum\nolimits_{i \in \mathcal{I}} \beta_{i,t}\hat{p}_{i,t}^w - E)\big] \\
%         & \leq B + V\mathbb{E}\big[\sum\nolimits_{i \in \mathcal{I}} \gamma_{i,t}(\check{p}_{i,t}^w -\hat{p}_{i,t}^w)\big]
%     \end{aligned}
% \end{equation}

Summing the above inequality over the time horizon $t \in \mathcal{T}$, dividing by $VT$ and driving $T$ to infinity yield the following:
\begin{equation}
    \begin{aligned}
       &\lim_{T \to \infty}\frac{1}{VT}\mathbb{E}[L(\mathbf{\Theta}_{T+1}^\ast)-L(\mathbf{\Theta}_{1}^\ast)] \\
       &+ \lim_{T \to \infty}\frac{1}{T}\sum_{t=1}^{T}\mathbb{E}[\gamma_{i,t}(\check{p}_{i,t}^\ast -\hat{p}_{i,t}^\ast)]\\
        & \leq \frac{B}{V}+ \lim_{T \to \infty}\frac{1}{T}\sum_{t=1}^{T}\mathbb{E}[\gamma_{i,t}(\check{p}_{i,t}^w -\hat{p}_{i,t}^w)]
    \end{aligned}
\end{equation}
where $B$ can be defined as:
\begin{align*}
        B &=  \frac{1}{2}\sum\nolimits_{j \in \mathcal{J}}\big[\max(A_{j}^2 , \sum_{i \in \Omega_j^f}M_{ji}^2) + \max ( \delta_{f,j}^2, \sum_{i \in \Omega_j^f}M_{ji}^2) \big] \\
         & + \frac{1}{2}\sum\nolimits_{i \in \mathcal{I}}\big[\max(\sum_{j \in \Omega_i^b}M_{ji}^2 , \check{P}_{i}^2) +  \max( \delta_{b,i}^2, \check{P}_{i}^2) \big]\\
         & + \sum\nolimits_{j \in \mathcal{J}}\big[ Z_{f,j}^{max}\delta_{f,j} \big] + \sum_{i \in \mathcal{I}}\big[Z_{b,i}^{max}\delta_{b,i} \big] \stepcounter{equation}\tag{\theequation}\\
         & +\frac{1}{2} \alpha \max( \sum\nolimits_{i \in \mathcal{I}} \bar{\beta}_{i}^2\hat{P}_{i}^2, E^2),
\end{align*}
and $A_j$, $\check{P}_i$, $\hat{P}_i$ and $\bar{\beta}_i$ are the upper bounds of $a_{j,t}$, $\check{p}_{i,t}$, $\hat{p}_{i,t}$ and ${\beta}_{i,t}$, respectively. 
Thus, we have:
\begin{equation}
    \Tilde{f} \leq \frac{B}{V}+ f^\ast
\end{equation}
Then we obtain:
\begin{equation}
    f_1^\ast - \Tilde{\varepsilon} \leq \frac{B}{V} + f^\ast
\end{equation}
Let $ \Tilde{\varepsilon} \to 0$, we have the following:
\begin{equation}
     f_1^\ast - f^\ast\leq \frac{B}{V}.
\end{equation}
This completes the proof.

\section{Formulation of baselines}\label{app_C}
SMA is implemented at time $t$ by solving the following:
\begin{subequations}
    \begin{align}
  & \max \;  \sum\nolimits_{i \in \mathcal{I}} \gamma_{i,t}(\Hat{p}_{i,t} -\Check{p}_{i,t}), \label{eq:P3_a}\\
&\textrm{s.t.} \quad \eqref{eq:lower_link_bandwidth} - \eqref{eq:lower_bound_limit}, \eqref{eq:upper_link_bandwidth} - \eqref{eq:upper_bound_limit},\\
& Q_{f,j,t} + a_{j,t} - \sum\nolimits_{i \in \Omega_{j}^f} \check{m}_{ji,t} \leq \overline{Q}_{f,j},\\
& Q_{b,i,t} + \sum\nolimits_{j \in \Omega_{i}^b} \check{m}_{ji,t} - \check{p}_{i,t} \leq \overline{Q}_{b,i},\\
&\sum\nolimits_{i \in \mathcal{I}} \beta_{i,t}\hat{p}_{i,t} \leq E,
    \end{align}
\end{subequations}
where $\overline{Q}_{f,j}$ and $\overline{Q}_{b,i}$ are the two associated upper limits.

AMA is implemented at time $t$ by solving the following:
\begin{subequations}
    \begin{align}
  & \max \;  \sum\nolimits_{i \in \mathcal{I}} \gamma_{i,t}(\Hat{p}_{i,t} -\Check{p}_{i,t}), \label{eq:P3_a}\\
&\textrm{s.t.} \quad \eqref{eq:lower_link_bandwidth} - \eqref{eq:lower_bound_limit}, \eqref{eq:upper_link_bandwidth} - \eqref{eq:upper_bound_limit},\\
& Q_{f,j,t} + a_{j,t} - \sum\nolimits_{i \in \Omega_{j}^f} \check{m}_{ji,t} \leq \overline{Q}_{f,j},\\
& Q_{b,i,t} + \sum\nolimits_{j \in \Omega_{i}^b} \check{m}_{ji,t} - \check{p}_{i,t} \leq \overline{Q}_{b,i},\\
&\sum\nolimits_{i \in \mathcal{I}} \beta_{i,t}\hat{p}_{i,t} \leq (E\times T - H_t)/(T-t+1),
    \end{align}
\end{subequations}
where $H_t$ is the accumulated carbon emissions up to current time $t$ and can be calculated by
\begin{equation}
    H_t = \begin{cases}
    0, &  t = 1 \\
    \sum_{\tau = 1}^{t-1} \sum_{ i \in \mathcal{I}}\beta_{i,\tau}{p}_{i,\tau}^d, &  t \ge 2.
\end{cases}
\end{equation}

\end{appendices}

\bibliographystyle{IEEEtran}
\bibliography{references}

\end{document}